\documentclass[12pt,leqno]{amsart}
\usepackage{amssymb}
\usepackage{amsmath,amssymb,color}
\numberwithin{equation}{section}
\newcommand{\DDD}{\mathcal{D}}
\newcommand{\DDDg}{(-A)^{\gamma}}

\newcommand{\ep}{\varepsilon}
\newcommand{\la}{\lambda}
\newcommand{\va}{\varphi}
\newcommand{\ppp}{\partial}

\newcommand{\pppb}{\partial_t^{\beta}}
\newcommand{\whwh}{\widehat}
\newcommand{\rrrr}{\longrightarrow}
\newcommand{\ddda}{d_t^{\alpha}}
\newcommand{\AAA}{\mathcal{A}}
\newcommand{\ABAB}{(-A)^{\beta}}
\newcommand{\AEAE}{(-A)^{\ep}}

\newcommand{\sumk}{\sum_{k=0}^\infty}

\newcommand{\pppa}{\partial_t^{\alpha}}

\newcommand{\R}{\mathbb{R}}
\newcommand{\C}{\mathbb{C}}
\newcommand{\N}{\mathbb{N}}

\newcommand{\www}{\widetilde}

\newcommand{\ooo}{\overline}
\newcommand{\OOO}{\Omega}

\newcommand{\PPPPO}{P_{\lambda}^{(p)}}
\newcommand{\PPPPT}{\widetilde{P}_{\lambda}^{(p)}}
\newcommand{\AAAAT}{\widetilde{A}_p}

\allowdisplaybreaks

\begin{document}
\title
[]
{
Initial boundary value problems for time-fractional evolution  
equations in Banach spaces 
}


\pagestyle{myheadings}

\author{
$^1$ Giuseppe Floridia, $^2$ Fikret G\"olgeleyen $^3$ Masahiro Yamamoto}
\thanks{
$^1$ 
Sapienza Universit\`a di Roma, 
Dipartimento di Scienze di Base e Applicate per l'Ingegneria
via Antonio Scarpa 16, 00161 Roma Italy
e-mail: {\tt giuseppe.floridia@uniroma1.it}
\\
$^2$ 
Department of Mathematics, Faculty of Science, 
Zonguldak B\"ulent Ecevit University, Zonguldak 67100 T\"urkiye
e-mail: {\tt f.golgeleyen@beun.edu.tr}
\\
$^3$ Graduate School of Mathematical Sciences, The University
of Tokyo, Komaba, Meguro, Tokyo 153-8914 Japan 
\\
Department of Mathematics, Faculty of Arts and Sciences, 
Zonguldak B\"ulent Ecevit University, Zonguldak 67100 Turkey
e-mail: {\tt myama@ms.u-tokyo.ac.jp}
}

\date{}

\begin{abstract}
We consider an initial value problem for time-fractional evolution equation
in Banach space $X$:
$$
\pppa (u(t)-a) = Au(t) + F(t), \quad 0<t<T.          \eqno{(*)}
$$
Here $u: (0,T) \rrrr X$ is an $X$-valued function defined in $(0,T)$,
and $a \in X$ is an initial value.
The operator $A$ satisfies a decay condition of resolvent which is 
common as a generator of analytic semigroup, and in particular,
we can treat a case $X=L^p(\OOO)$ over a bounded domain $\OOO$
and a uniform elliptic operator $A$ within our framework.

First we construct a solution operator $(a, F) \rrrr u$ by means of 
$X$-valued Laplace transform, and we establish the well-posedness of (*)
in classes such as weak solution and strong solutions.  
We discuss also mild solutions local in time for semilinear 
time-fractional evolution equations.
Finally we apply the result on the well-posedness to an inverse problem of 
determining an initial value 
and we establish the uniqueness for the inverse problem.
\end{abstract} 
\baselineskip 18pt

\maketitle
\section{Introduction and formulation}

Let $X$ be a Banach space over $\C$ and let $A: \DDD(A) \longrightarrow
X$ be a densely defined closed linear operator satisfying 
{\it Condition ($\AAA$)} stated later.  

For $u=u(t): [0,T] \rrrr X$, we introduce an initial value problem for
a time-fractional evolution equation in $X$:
$$
\left\{ \begin{array}{rl}
& \ddda u(t) = Au(t) + F(t) \quad \mbox{for $t>0$ in $X$,}\\
& u(0) = a.
\end{array}\right.
                                          \eqno{(1.1)}
$$
Here we set $u'(t):= \frac{du}{dt}(t)$ and for $0<\alpha<1$, 
the Caputo derivative $\ddda u$ is defined by 
$$
\ddda u(t) = \frac{1}{\Gamma(1-\alpha)}\int^t_0 (t-s)^{-\alpha}
u'(s) ds,
$$
provided that the right-hand side can be defined (e.g., $u \in C^1([0,T];X)$).

The main purpose of this article is to prove the well-posedness of (1.1), i.e.,
the unique existence of a solution $u$ to (1.1) with an estimate by 
suitable norms of $a$ and $F$.  The initial value problem (1.1) in 
a Banach space, describes also 
initial boundary value problems for time-fractional diffusion equations by 
adequate choices of $A$ (an elliptic operator, for example)
The well-posedness of the initial boundary value problem is the primary 
step for studies of time-fractional diffusion equations and here we mainly 
discuss in terms of time-fractional evolution equations (1.1).

There are very many works on the well-posedness in the case where $X$ is 
a Hilbert space and $A$ admits an orthonormal basis composed of 
the eigenvectors of $A$ (for example, $A$ is self-adjoint with compact 
resolvent), and we can refer to Jin \cite{J}, Kubica and Yamamoto \cite{KY},
Kubica, Ryszewska and Yamamoto \cite{KRY}, Zacher \cite{Za}.  
We further refer to Achache \cite{Ach1}, \cite{Ach2}.
Here we do not
intend any comprehensive references and see the references therein 
as for further literature.

On the other hand, for general Banach space $X$, there are few results 
which admit feasible theories for various applications such as 
nonlinear equations, inverse problems.
As works focusing on well-posedness, see for example, 
Pr\"uss \cite{Pr}, Zacher \cite{Za}.
The monograph \cite{Pr} treats general equations including (1.1), but 
when one restrict equations to more concrete equations such as 
(1.1), it is not necessarily clear how much one can sharpen the results.
As more recent works, see Henr\'iquez, Mesquita ad Pozo \cite{HMP},
Kexue, Jigen and Junxiong \cite{KJJ}, and also here we do not aim at 
any comprehensive references.

In this article, we will establish the well-posedness of 
initial value problems for time-fractional evolution equations in 
Banach space $X$ and $L^q$-space in time $t$, where $1\le q \le \infty$.
Our method is a modification of the classical method for
construction of the analytic semigroups
via the Laplace transform.  As monographs, we can refer for example to 
Kato \cite{Ka}, Pazy \cite{Pa}, Tanabe \cite{Ta}, Yagi \cite{Ya}.
Such a construction is very classical for $\alpha=1$, but there are no 
existing works for $0<\alpha<1$, to the best knowledge of the authors.

The construction is feasible for several applications and so
we will study some of the applications.

We now formulate the initial value problem.
First we need to define fractional derivatives in 
the Lebesgue spaces and introduce a class of operators $A$ and function spaces.
Throughout this article, we assume that the indexes of the function spaces
satisfy 
$$
1\le p, q \le \infty,
$$
if we do not specify.
We define $\pppa$ in $L^q(0,T;X)$ where $X$ is a Banach space $X$.

For $\beta > 0$, we set 
$$
J^{\beta}v(t):= \frac{1}{\Gamma(\beta)}\int^t_0 (t-s)^{\beta-1}
v(s) ds, \quad v\in L^q(0,T;X).            \eqno{(1.2)}
$$
Then, we can easily prove that
$J^{\beta}: L^q(0,T;X) \longrightarrow J^{\beta}L^q(0,T;X)$ is injective and 
surjective.  We define a function space  
$$
W_{\beta,q}(0,T;X) := J^{\beta}L^q(0,T;X)                \eqno{(1.3)}
$$
with the norm $\Vert v\Vert_{W_{\beta,q}(0,T;X)}
:= \Vert (J^{\beta})^{-1}v \Vert_{L^q(0,T;X)}$.
We can verify that $W_{\beta,q}(0,T;X)$ is a Banach space.

We define a time-fractional derivative in $W_{\beta,q}(0,T;X)$ by 
$$
\ppp_t^{\beta} := (J^{\beta})^{-1}, \quad
\DDD(\ppp_t^{\beta}) = J^{\beta}L^q(0,T;X).
                                                   \eqno{(1.4)}
$$
Henceforth we denote $J^{-\beta}:= (J^{\beta})^{-1}$ not only 
algebraically but also including the topology.

Then 
\\
{\bf Lemma 1.1.}
\\
{\it
$\pppb : \DDD(\pppb) = W_{\beta,q}(0,T;X)\, \longrightarrow
\, L^q(0,T;X)$ is a closed operator.
}
\\
{\bf Proof.}
Let $u_n \in \DDD(\pppb)$ for $n\in \N$ and $u_n \longrightarrow u$ in 
$L^q(0,T;X)$, $\pppb u_n \longrightarrow w$ in $L^q(0,T;X)$ with some 
$u, w \in L^q(0,T;X)$.
Since $u_n\in \DDD(\pppb)$, we can find $w_n \in L^q(0,T;X)$ such that 
$u_n = J^{\beta}w_n$, that is, $w_n = \pppb u_n$,
Then $w_n \longrightarrow w$ and $J^{\beta} w_n \longrightarrow u$
in $L^q(0,T;X)$.
Since $J^{\beta}: L^q(0,T;X) \longrightarrow L^q(0,T;X)$ is bounded 
by the Young inequality on the convolution, we know that 
$J^{\beta} w_n \longrightarrow J^{\beta}w = u$ in $L^q(0,T;X)$.
This means that $w = \pppb u$ and $u \in \DDD(\pppb)$.
The proof of Lemma 1.1 is completed.
$\blacksquare$
\\

Henceforth by $\rho(A)$ we denote the resolvent set of an operator 
$A$ and 
$$
\Sigma_{\gamma}:= \{ z\in\C;\, \vert \mbox{arg}\, z\vert < \gamma,
\,\, z\ne 0\}
$$
with $0 < \gamma < \pi$.
 
For the operator $A$, we pose 
\\
{\bf Condition ($\AAA$):}
\\
(i) $A$ is a closed linear operator in $X$ such that $\DDD(A)$ is dense
in $X$.
\\
(ii) There exists a constant $\gamma \in \left( \frac{\pi}{2}, \, \pi\right)$
such that $\Sigma_{\gamma} \subset \rho(A)$.
\\
(iii) For arbitrarily given $\ep \in (0, \gamma)$, there exists a constant 
$C = C_{\ep}>0$ such that 
$$
\Vert (\la - A)^{-1}\Vert \le \frac{C_{\ep}}{\vert \la\vert}
\quad \mbox{for all $\la \in \Sigma_{\gamma-\ep}$}.             \eqno{(1.5)}
$$
\\
(iv) $0 \in \rho(A)$.
\\

In this article, we consider the case $0<\alpha<1$.

Under these frameworks, we formulate the initial value problem for
time-fractional evolution equation:
$$
\pppa (u(t)-a) = Au(t) + F(t) \quad \mbox{for $t>0$ in $X$.}
                                                            \eqno{(1.6)}
$$
We aim at the well-posedness of (1.6) in a Banach space $X$, which is not 
necessarily neither a Hilbert space nor a reflexive space.
Our main methodology is based on construction of solution mapping:
$a \mapsto u(t)$ with $F=0$ and the theory for 
$X$-valued Laplace transforms (e.g., Arendt, Batty, Hieber and
Neubrander \cite{ABHN}).   We note that our solution mapping 
for $0<\alpha<1$ corresponds to the semigroup $e^{tA}$ for the case
$\alpha=1$, although our solution mapping has no properties as
semigroups.

The initial value problem (1.6) includes a classical initial boundary 
value problem for a time-fractional diffusion equation.  More precisely,
let $\OOO \subset \R^d$ be a bounded domain with smooth boundary $\ppp\OOO$.
We formally define an elliptic operator $\mathcal{A}$ by
$$
\mathcal{A}v(x) := \sum_{k,\ell=1}^d \ppp_k(a_{k\ell}(x)\ppp_{\ell}v(x)) 
+ \sum_{k=1}^d b_k(x)\ppp_kv(x) + c(x)v(x),                        \eqno{(1.7)}
$$
where $a_{k \ell} = a_{\ell k} \in C^2(\ooo{\OOO})$, $b_k$, 
$c \in C(\ooo{\OOO})$, $c \le 0$ in $\OOO$, and we assume that 
there exists a constant $\kappa>0$ such that 
$$
\sum_{k,\ell=1}^d a_{k\ell}(x)\xi_k\xi_{\ell} 
\ge \kappa \sum_{k=1}^d \xi_k^2, \quad 
x\in \ooo{\OOO}, \, \xi_1, ..., \xi_d \in \R.    \eqno{(1.8)}
$$
In order to set up $\mathcal{A}$ within the framework of the
operator theory, to $\mathcal{A}$ we attach the domain 
$$
\DDD(A) = \{ u \in W^{2,p}(\OOO);\, 
u\vert_{\ppp\OOO} = 0\}
$$
and by $A$ we define such an elliptic operator with the domain 
$\DDD(A)$.

Then, the condition $c\le 0$ implies that $0 \in \rho(A)$, and 
for $1<p<\infty$, by means of Theorem 3.2 (p.213) in \cite{Pa} for example, 
it is known that
$A$ satisfies {\it Condition ($\AAA$)}.
\\

This article is composed of eight sections.
In Section 2, we construct a solution operator $G(t): [0,T] \, \rrrr X$ for
$\pppa (u(t) - a) = Au(t)$ for $t>0$ to establish the well-posedness of
(1.6) with $F=0$.  In Section 3, we consider 
$\pppa u(t) = Au(t) + F(t)$ for $t>0$.
Section 4 complies the results in Sections 2 and 3 and state 
the main result on the well-posedness for (1.6) for $(a,F)$ including 
in a suitable space.
Then, for demonstrating the feasibility of our results, we discuss:
\begin{itemize}
\item
Sectoin 5. Time local existence of solutions to semilinear equations
\item
Section 6. Weak and smoother solutions
\item
Section 7. An application to inverse problem of determining initial values 
\item
Section 8. Concluding remarks
\end{itemize}

\section{construction of a solution operator in the case $F=0$}

We consider
$$
\pppa (u(t) - a) = Au(t) \quad \mbox{for $t>0$ in $X$.}    \eqno{(2.1)}
$$
Henceforth we define the Laplace transform of $u: [0,\infty) \rrrr X$ by 
$$
(Lu)(\la):= \int^{\infty}_0 e^{-\la t} u(t) dt,
$$
provided that the right-hand side exists for $\la \in \C$.
We recall that $\gamma \in \left( \frac{\pi}{2}, \, \pi\right)$ is 
specified in {\it Condition ($\AAA$)}.
\\

We set
$$
G(z)a:= \frac{1}{2\pi i}\int_{\Gamma} 
e^{\la z}\la^{\alpha-1}(\la^{\alpha}-A)^{-1}a d\la \quad 
\mbox{for $a \in X$ and $z\in \Sigma_{\gamma-\frac{\pi}{2}}$}.
                               \eqno{(2.2)}
$$
Here we can choose a path $\Gamma\subset \rho(A)$ as follows:
$\Gamma$ starts at $+\infty e^{-i\gamma}$ to reach near 
the origin $0\in \C$ and then avoiding $0$, surrounds $0$, and again
goes to $+\infty e^{i\gamma}$.
In surrounding $0$, the path remains so that $\{ \vert \mbox{Re}\, z\vert;\, 
z\in \Gamma\}$ is sufficiently small which implies  $\Gamma \subset \rho(A)$.  
Since $\Gamma$ does not intersect $\{ z\le 0\}$, we see that $\la^{\alpha-1}
(\la^{\alpha} - A)^{-1}a$ is holomorphic in a neighborhood of 
$\Gamma$.  Hence, by means of Cauchy's integral theorem, we can verify that 
$G(z)a$ is invariant under transformations of $\Gamma$ satisfying the 
above conditions.  For example, as $\Gamma$ we can choose the following 
two choices which provide the same value $G(z)a$:
\\
(a) $z$-dependent $\Gamma:= \Gamma_1 \cup \Gamma_2 \cup \Gamma_3$.
$$
\Gamma_1 = \Gamma_1(z) := \left\{\rho e^{-i\gamma};\, \rho 
> \frac{1}{\vert z\vert} \right\},
\quad
\Gamma_2:= \left\{\frac{1}{\vert z\vert}e^{i\theta};\, 
-\gamma \le \theta \le \gamma
\right\}, 
$$
$$
\Gamma_3:= \left\{\rho e^{i\gamma};\, \rho > \frac{1}{\vert z\vert} \right\}.
                                                     \eqno{(2.3)}
$$
\\
(b) $z$-independent $\Gamma:= \Gamma^1 \cup \Gamma^2 \cup \Gamma^3$.
$$
\Gamma^1 = \{\rho e^{-i\gamma};\, \rho > \ep\},
\quad
\Gamma^2:= \{ \ep e^{i\theta};\, -\gamma \le \theta \le \gamma\},
\quad
\Gamma^3:= \{\rho e^{i\gamma};\, \rho > \ep\},
                                                     \eqno{(2.4)}
$$
where the constant $\ep>0$ is sufficiently small.

In this article, we mainly choose $\Gamma$ defined by (2.3).
\\

In this section, we firstly establish several estimates of $G(t)a$ and 
related equations, which guarantee the existence of $G(t)a$ in suitable spaces.
In terms of {\it Condition ($\AAA$)}, we can define a fractional power 
$(-A)^{\beta}$ of $-A$ (e.g., \cite{Pa}, \cite{Ta}).
The following estimates are proved.
\\
{\bf Lemma 2.1.}
\\
{\it
Let $0 \le \beta \le 1$.  Then the operator $(-A)^{\beta}$ is closed in 
$X$, and there exists a constant $C_{\beta} > 0$ 
such that 
$$
\Vert (-A)^{\beta}(A-\la)^{-1}a\Vert
\le C_{\beta}\vert \la\vert^{\beta-1}\Vert a\Vert
$$
if $\la \in \Sigma_{\gamma}$.
}
\\
\vspace{0.2cm}
\\
The proof is found e.g., Corollary (p.39) in \cite{Ta}.
\\
\vspace{0.2cm}

We state the estimation of $G(t)a$ which we will prove in this section.
\\
{\bf Proposition 2.1.}
\\
{\it 
Let $0 \le \beta \le 1$.
\\
(i) $AG(z)a = G(z)Aa$ for all $z \in \Sigma_{\gamma-\frac{\pi}{2}}$ and 
$a \in \DDD(A)$.
\\
(ii) There exists a constant $C = C(\beta) > 0$ such that 
$$
\Vert (-A)^{\beta}G(z)a\Vert \le C\vert z\vert^{-\alpha\beta}\Vert a\Vert
$$
for all $z \in \Sigma_{\gamma-\frac{\pi}{2}}$ and $a \in X$.
\\
(iii) $G(z)a$ is holomorphic in $z\in \Sigma_{\gamma-\frac{\pi}{2}}$.
}
\\
{\bf Proof.}
\\
(i) In terms of the closedness of the operator $A$, we can verify
\\
{\bf Lemma 2.2.}
\\
{\it
For an $X$-valued function $V: \Gamma \rrrr X$, we assume that 
$V(\la) \in \DDD(A)$ for each $\la \in \Gamma$ and
$AV(\cdot) \in C(\Gamma)$, $\Vert AV(\cdot)\Vert \in L^1(\Gamma)$.
Then 
$$
\int_{\Gamma} V(\la) d\la \in \DDD(A), \quad
A\left( \int_{\Gamma} V(\la) d\la\right)
= \int_{\Gamma} AV(\la) d\la.
$$
}
\\

We apply Lemma 2.2 to $V(\la) := \frac{1}{2\pi i}e^{\la z}
\la^{\alpha-1}(\la^{\alpha} - A)^{-1}a$.
We note that $\Gamma$ is defined by (2.3).  Since 
$$
AV(\la)
:= \frac{1}{2\pi i} e^{\la z}\la^{\alpha-1} (\la^{\alpha}-A)^{-1}Aa,
$$ 
by (1.5) with $\ep=0$, we obtain
$$
\Vert A^jV(\la)\Vert \le \frac{1}{2\pi}\sup_{\la\in \Gamma}
\vert e^{\mbox{Re}\,(\la z)}\vert \vert \la\vert^{\alpha-1}
\frac{\Vert A^ja\Vert}{\vert \la\vert^{\alpha}}
\le C \sup_{\la\in \Gamma}
\vert e^{\mbox{Re}\,(\la z)}\vert \frac{\Vert A^ja\Vert}{\vert \la\vert},
\quad \la \in \Sigma_{\gamma}
$$
for $j=0,1$.

We estimate Re $(\la z)$ for $z \in \Sigma_{\gamma-\frac{\pi}{2}}$.  Then,
$$
z = \vert z\vert e^{i\psi}, \quad 
 -\gamma+\frac{\pi}{2} < \psi < \gamma-\frac{\pi}{2}.     \eqno{(2.5)}
$$
On $\Gamma_1 \cup \Gamma_3$, we have 
$\la = \vert \la\vert e^{\pm i\gamma} 
= \vert \la\vert(\cos \gamma \pm i \sin \gamma)$ 
and so
$$
\mbox{Re}\, (\la z) = \mbox{Re}\, (\vert \la\vert
\vert z\vert e^{i\psi}e^{\pm i\gamma})
= \vert \la\vert \vert z\vert \cos (\gamma \pm \psi).
$$
For $-\gamma+\frac{\pi}{2} < \psi < \gamma-\frac{\pi}{2}$, we can directly 
verify that $\frac{\pi}{2} < \gamma \pm \psi < \frac{3}{2}\pi$.
Therefore, we can choose a constant $\delta>0$ such that 
$\cos (\gamma \pm \psi) \le -\delta$ for all 
$\vert \psi \vert < \gamma-\frac{\pi}{2}$.

For $\la \in \Gamma_2$, we see that 
$\la = \frac{1}{\vert z\vert}(\cos \theta + i \sin \theta)$ for 
$-\gamma < \theta < \gamma$, and so 
$\la z = (\cos \theta + i\sin\theta)(\cos \psi + i\sin\psi)$, 
which implies
$$
\mbox{Re}\, (\la z) =  \cos\theta \cos\psi - \sin \theta \sin \psi
= \cos (\theta - \psi).
$$
Thus
$$
\left\{ \begin{array}{rl}
& \mbox{Re}\, (\la z) \le -\vert \la\vert \vert z\vert \delta \quad 
\mbox{for all $\la \in \Gamma_1 \cup \Gamma_3$ and 
$z \in \Sigma_{\gamma-\frac{\pi}{2}}$}, \\
& \mbox{Re}\, (\la z) = \cos (\theta - \psi) \quad 
\mbox{for all $\la \in \Gamma_2$ and 
$z \in \Sigma_{\gamma-\frac{\pi}{2}}$}.
\end{array}\right.
                         \eqno{(2.6)}
$$
In particular, we have
$$
\left\{ \begin{array}{rl}
& \mbox{Re}\, (\la t) \le -\vert \la\vert t \delta \quad 
\mbox{for all $\la \in \Gamma_1 \cup \Gamma_3$ and $t > 0$}, \\
& \mbox{Re}\, (\la t) 
= \cos \theta \quad \mbox{for all $\la \in \Gamma_2$ and $t>0$}.
\end{array}\right.
                             \eqno{(2.7)}
$$
Therefore, for $t>0$, we obtain
$$
\Vert A^jV(\la)\Vert \le
\left\{ \begin{array}{rl}
& Ce^{-\vert \la\vert t\delta}\vert z\vert \Vert A^ja\Vert,\quad
\la \in \Gamma_1 \cup\Gamma_3,  \\
& Ce \vert z\vert \Vert A^ja\Vert,\quad 
\la \in \Gamma_2
\end{array}\right.
$$
and for any fixed $t>0$, we have $\Vert A^jV(\cdot)\Vert 
\in L^1(\Gamma)$.  The rest assumptions of Lemma 2.2 are readily 
verified, so that the proof of Proposition 2.1 (i) is completed.
$\blacksquare$
\\
{\bf Proof of Proposition 2.1 (ii)}
\\
By (2.2) and the closedness of the operator $(-A)^{\beta}$, we can justify 
$$
(-A)^{\beta}G(z)a = \frac{1}{2\pi i}\int_{\Gamma} 
e^{\la z}\la^{\alpha-1}(-A)^{\beta}(\la^{\alpha}-A)^{-1}a d\la, \quad 
z\in \Sigma_{\gamma-\frac{\pi}{2}},
$$
and so (2.7) yields
\begin{align*}
& \Vert (-A)^{\beta}G(z)a \Vert 
\le C \int_{\Gamma} 
\vert e^{\la z}\vert \vert \la\vert^{\alpha-1}
\Vert (-A)^{\beta}(\la^{\alpha}-A)^{-1}a\Vert \vert d \la\vert  \\
=& C\left( \int_{\Gamma_1} + \int_{\Gamma_2} + \int_{\Gamma_3}\right) 
\vert e^{\mbox{Re}\, (\la z)}\vert \vert \la\vert^{\alpha-1}
\Vert (-A)^{\beta}(\la^{\alpha}-A)^{-1}a\Vert \vert d\la\vert\\
\le& C\int_{\Gamma_1} e^{-\vert \la\vert\vert z\vert \delta}
\vert \la\vert^{\alpha-1} \Vert (-A)^{\beta}(\la^{\alpha}-A)^{-1}a\Vert
\vert d\la\vert 
+ C\int_{\Gamma_2} e^{\cos(\theta-\psi)}
\vert \la\vert^{\alpha-1} \Vert (-A)^{\beta}(\la^{\alpha}-A)^{-1}a\Vert
\vert d\la\vert \\
+& C\int_{\Gamma_3} e^{-\vert \la\vert \vert z\vert \delta}
\vert \la\vert^{\alpha-1} \Vert (-A)^{\beta}(\la^{\alpha}-A)^{-1}a\Vert
\vert d\la\vert \\
=: &I_1(z) + I_2(z) + I_3(z).
\end{align*}
\\
{\bf Estimation of $I_1(z)$.}
We set $\la = \rho e^{-i\gamma}$ with $\frac{1}{\vert z\vert} < \rho < \infty$.
Then $d\la = e^{-i\gamma}d\rho$ and $\vert d\la\vert = d\rho$.
Using Lemma 2.1, we obtain
$$
I_1(z) 
= C\int^{\infty}_{\frac{1}{\vert z\vert}} e^{-\rho \vert z\vert\delta}
\rho^{\alpha-1}C\rho^{\alpha(\beta-1)} \Vert a\Vert \, d\rho
= C\int^{\infty}_{\frac{1}{\vert z\vert}} e^{-\rho\vert z\vert \delta}
\rho^{\alpha\beta-1} d\rho \Vert a\Vert.
$$
Setting $\eta = \rho \vert z\vert$, we have $d\eta = \vert z\vert d\rho$, 
and so 
$$
 \int^{\infty}_{\frac{1}{\vert z\vert}} e^{-\rho \vert z\vert \delta}
\rho^{\alpha\beta-1} d\rho
= \int^{\infty}_1 e^{-\delta \eta}\left( \frac{\eta}{\vert z\vert}\right)
^{\alpha\beta-1} 
\frac{1}{\vert z\vert} d\eta 
= \vert z\vert^{-\alpha\beta}\int^{\infty}_1 e^{-\delta\eta} 
\eta^{\alpha\beta-1} d\eta
=: C\vert z\vert^{-\alpha\beta}.
$$
Consequently,
$$
\vert I_1(z)\vert \le C\vert z\vert^{-\alpha\beta}\Vert a\Vert.
$$
\\
{\bf Estimation of $I_2(z)$.}
We set $\la = \frac{1}{\vert z\vert}e^{i\theta}$ where $\theta: - \gamma 
\longrightarrow \gamma$.  Then $d\la = \frac{1}{\vert z\vert}
ie^{i\theta}d\theta$,
$\vert \la \vert = \frac{1}{\vert z\vert}$ and $\vert d\la\vert 
= \frac{1}{\vert z\vert}d\theta$.
Hence,
\begin{align*}
& I_2(z) = C\int^{\gamma}_{-\gamma} 
e \left( \frac{1}{\vert z\vert}\right)^{\alpha-1}
C \left( \frac{1}{\vert z\vert}\right)^{\alpha(\beta-1)}
\frac{1}{\vert z\vert} d\theta \Vert a\Vert \\
=& \frac{Ce}{\vert z\vert^{\alpha\beta}} \int^{\gamma}_{-\gamma}
 d\theta \Vert a\Vert
\le C\vert z\vert^{-\alpha\beta} \Vert a\Vert
= \frac{2C\gamma e}{\vert z\vert^{\alpha\beta}}\Vert a\Vert.
\end{align*}
For $I_3(z)$, we can estimate similarly.
Thus, the proof of Proposition 2.1 (ii) is complete.

Part (iii) follows from Theorem 2.6.1 in \cite{ABHN}, which asserts the 
equivalence between the holomorphy of $G(z)a$ and the Laplace transform
$L(G(z)a)(\la)$.
Thus the proof of Proposition 2.1 is complete.
$\blacksquare$
\\

Now we state the main result in Section 2.
\\
{\bf Theorem 2.1}
\\
{\it
We arbitrarily fix a constant $\mu$ such that 
$$
1 - \frac{1}{q\alpha} < \mu < 1.        \eqno{(2.8)}
$$
Then
$$
Ga - a \in W_{\alpha,q}(0,T;X) \cap L^{\infty}(0,\infty;X) \quad
\mbox{for all $T>0$,}
                                                           \eqno{(2.9)}
$$
and
$$
\pppa (G(t)a-a) = AG(t)a \quad \mbox{for $a\in \DDD((-A)^{\mu})$} \quad
\mbox{for all $0 < t < T$}.
                                         \eqno{(2.10)}
$$
}
We note that (2.9) means $\pppa (Ga - a) \in L^q(0,T;X)$.
\\
{\bf Proof of Theorem 2.1}
\\
{\bf First Step}
\\
First we show
\\
\\
{\bf Lemma 2.3}
\\
{\it
Let $u \in W_{\alpha,q}(0,T;X)$ for any $T>0$ and $\Vert u(\cdot,t)\Vert
= O(t^m)$ as $t \to \infty$ with some $m\ge 0$.
Then, $(L\pppa u)(\la)$ exists for $\la > 0$ and
$$
(L\pppa u)(\la) = \la^{\alpha}(Lu)(\la) \quad \mbox{for $\la > 0$}.
$$
}
\\
{\bf Proof of Lemma 2.3}
\\
We have
$$
\pppa u(t) = \frac{1}{\Gamma(1-\alpha)}\frac{d}{dt}
\int^t_0 (t-s)^{-\alpha}u(s) ds,
$$
and then
\begin{align*}
& \int^T_0 \frac{d}{dt}\left( 
\int^t_0 (t-s)^{-\alpha}u(s) ds\right) e^{-\la t} dt\\
=& \left[ \left(\int^t_0 (t-s)^{-\alpha}u(s) ds\right) e^{-\la t}
\right]^{t=T}_{t=0}
+ \int^T_0 \left( \int^t_0 (t-s)^{-\alpha}u(s) ds\right)
\la e^{-\la t} dt.
\end{align*}
Then
$$
\lim_{T\to \infty} \left(\int^T_0 (T-s)^{-\alpha}u(s) ds\right) e^{-\la T} = 0.
$$
Indeed, 
\begin{align*}
& \left\Vert \int^T_0 (T-s)^{-\alpha}u(s) ds \right\Vert e^{-\la T} 
\le \left( \int^T_0 (T-s)^{-\alpha}\Vert u(s)\Vert ds \right)e^{-\la T} \\
\le & C\left( \int^T_0 (T-s)^{-\alpha}T^m ds \right)e^{-\la T} 
= CT^me^{-\la T}\frac{T^{1-\alpha}}{1-\alpha} \longrightarrow 0
\end{align*}
as $T \to \infty$.

Moreover, 
$$
 \la\int^T_0 \left(\int^t_0 (t-s)^{-\alpha}u(s) ds\right) e^{-\la t} dt
= \la\int^T_0 \left(\int^T_s (t-s)^{-\alpha}e^{-\la t} dt\right) u(s) ds.
$$
Here setting $\eta:= t-s$, that is, $t = \eta + s$, we have
$$
\int^T_s (t-s)^{-\alpha}e^{-\la t} dt
= e^{-\la s}\int^{T-s}_0 \eta^{-\alpha}e^{-\la\eta} d\eta
\longrightarrow e^{-\la s}\frac{\Gamma(1-\alpha)}{\la^{1-\alpha}}
$$
as $T \to \infty$.
Hence,
$$
\la\int^T_0 \left(\int^t_0 (t-s)^{-\alpha}u(s) ds\right) e^{-\la t} dt
\, \longrightarrow \, \la^{\alpha}\Gamma(1-\alpha)\whwh{u}(\la).
$$
Thus the proof of Lemma 2.3 is complete.
$\blacksquare$
\\

Now we can prove
$$
A(LGa)(\la) = L(AGa)(\la) \quad \mbox{for all $\la > 0$ and
$a\in \DDD(A)$.}                    \eqno{(2.11)}
$$
{\bf Verification of (2.11).}
\\
We arbitrarily fix $\la > 0$.  In terms of $a\in \DDD(A)$, Proposition 2.1 (i)
yields $G(t)a \in \DDD(A)$ and $AG(t)a = G(t)Aa$ for all $t>0$.
Moreover $A(G)a\in C([0,\infty);X) \cap L^{\infty}(0,\infty;X)$ and
$e^{-\la t}G(t)a \in \DDD(A)$ for each $t>0$.
Then $e^{-\la t}AG(t)a \in C_t([0,\infty);X) \cap L_t^1(0,\infty;X)$ for 
arbitrarily fixed $\la > 0$.  
Here $C_t([0,\infty);X)$ and $L_t^1(0,\infty;X)$ 
mean the corresponding function
spaces with the independent variable $t$.

A similar application of Lemma 2.2 yields  
\begin{align*}
& AL(Ga)(\la) = A\int^{\infty}_0 e^{-\la t}G(t)a dt
= \int^{\infty}_0 A(e^{-\la t}G(t)a) dt
= \int^{\infty}_0 e^{-\la t}(AG(t)a) dt\\
=& L(AGa)(\la) \quad \mbox{for all $\la > 0$.}
\end{align*}
Thus the verification of (2.11) is complete.
$\blacksquare$

{\it Condition ($\AAA$)} yields 
$$
\sup_{\la \in \Sigma_{\gamma}} \Vert \la \la^{\alpha-1}(\la^{\alpha} - A)^{-1}a
\Vert < \infty,
$$
and so we can apply Theorem 2.6.1 in \cite{ABHN} to conclude that 
$$
G(z)a = \frac{1}{2\pi i}\int_{\Gamma} e^{\la z} \la^{\alpha-1}
(\la^{\alpha} - A)^{-1} a d\la
$$
is holomorphic in $\Sigma_{\gamma-\frac{\pi}{2}}$ and
$$
(LGa)(\la) = \la^{\alpha-1}(\la^{\alpha}-A)^{-1}a \quad
\mbox{for all $\la > 0$}.        \eqno{(2.12)}
$$

By (2.12), we have
$$
 \la^{\alpha}L(Ga)(\la) - A(LGa)(\la)
= (\la^{\alpha}-A)L(Ga)(\la) = \la^{\alpha-1}a
$$
for $\la > 0$.  Therefore, (2.11) yields
$$
\la^{\alpha}L(Ga)(\la) - L(AGa)(\la) = \la^{\alpha-1}a
$$
for $\la > 0$.  Since $(La)(\la) = \la^{-1}a$ by direct calculations,
we obtain
$$
\la^{\alpha}L(Ga-a)(\la) = L(AGa)(\la), \quad \la > 0.
                                   \eqno{(2.13)}
$$
Setting $w(t) := AG(t)a$, we rewrite (2.13) as
$$
Lw(t) = \la^{\alpha}L(Ga-a)(\la), \quad \la > 0.
                                       \eqno{(2.14)}
$$
In view of Proposition 2.1 (i) and (ii), we have
$$
w(t) = G(t)Aa \in L^{\infty}(0,\infty;X) \subset L^q(0,T;X) \quad
\mbox{for any $T>0$.}
$$
Hence, setting 
$$
v:= J^{\alpha}w \in W_{\alpha,q}(0,T;X) \quad \mbox{for any $T>0$},
                                  \eqno{(2.15)}
$$
we see that $w = \pppa v$ and
\begin{align*}
& \Vert v(t)\Vert_X = \left\Vert \frac{1}{\Gamma(\alpha)}
\int^t_0 (t-s)^{\alpha-1} w(s) ds\right\Vert \\
\le& C\int^t_0 (t-s)^{\alpha-1}\Vert w(s)\Vert ds 
= C\Vert w\Vert_{L^{\infty}(0,\infty;X)}t^{\alpha} \quad
\mbox{for all $t>0$.}
\end{align*}
Therefore, Lemma 2.3 yields
$$
Lw(\la) = L(\pppa v)(\la) = \la^{\alpha}Lv(\la), \quad \la > 0.
                                                    \eqno{(2.16)}
$$
In terms of (2.14), we have
$$
\la^{\alpha}Lv(\la) = \la^{\alpha}L(Ga - a)(\la), \quad \la > 0,
$$
that is,
$$
(Lv)(\la) = L(Ga - a)(\la), \quad \la > 0.
$$
The injectivity of the Laplace transform yields
$$
v(t) = G(t)a - a, \quad t>0.            \eqno{(2.17)}
$$
By (2.15), we see that $v \in W_{\alpha,q}(0,T;X)$ for all $T>0$.
Thus we have proved that $Ga - a \in W_{\alpha,q}(0,T;X)$ in (2.9).

Finally (2.16) and (2.17) yield
\begin{align*}
& Lw(\la) = L(\pppa v)(\la)
= L(\pppa (Ga - a))(\la)
= \la^{\alpha}(Lv)(\la)\\
=& \la^{\alpha}L(Ga-a)(\la), \quad \la > 0.
\end{align*}
Hence, (2.13) implies
$$
L(\pppa (Ga-a))(\la) = L(AGa)(\la), \quad \la>0.
$$
Again the injectivity of the Laplace transform yields
$AG(t)a = \pppa (G(t)a - a)$ for $t>0$.
Thus the proof of (2.9) and (2.10) is completed for $a \in \DDD(A)$.
$\blacksquare$
\\
{\bf Second Step.}
\\
We will prove (2.9) and (2.10) for $a \in \DDD((-A)^{\mu})$, where
$\mu$ satisfies (2.8).
Since $\DDD(A)$ is dense in $X$, for each $a \in \DDD((-A)^{\mu})$, setting
$b:= (-A)^{\mu}a\in X$, we can find a sequence $b_n\in \DDD(A)$, $n\in \N$ 
such that $b_n \rrrr b$ in $X$ as $n \to \infty$.
Therefore, $(-A)^{-\mu}b_n \rrrr (-A)^{-\mu}b = a$ in $\DDD((-A)^{\mu})$
as $n\to \infty$.
By $b_n \in \DDD(A)$, we have 
$A((-A)^{-\mu}b_n) = (-A)^{-\mu}(Ab_n) \in X$, that is,
$a_n:= (-A)^{-\mu}b_n \in \DDD(A)$.  Hence, or each $a \in \DDD((-A)^{\mu})$,
we can find a sequence $a_n \in \DDD(A)$, $n\in \N$ such that 
$(-A)^{\mu}a_n \rrrr (-A)^{\mu}a$ in $X$ as $n\to \infty$, that is,
$$
a_n \rrrr a \quad \mbox{in $\DDD((-A)^{\mu})$ as $n\to \infty$}.
                                         \eqno{(2.18)}
$$

As is already proved, since $a_n \in \DDD(A)$, 
we see that $G(t)a_n - a_n$ satisfies 
(2.9) and (2.10).  By (2.18) and Proposition 2.1 (ii), we can obtain
\begin{align*}
& \Vert A(G(t)(a_n-a))\Vert = \Vert (-A)^{1-\mu}G(t)(-A)^{\mu}(a_n-a)\Vert\\
\le & Ct^{-\alpha(1-\mu)}\Vert (-A)^{\mu}(a_n-a)\Vert.
\end{align*}
Consequently,
\begin{align*}
& \Vert AG(a_n-a)\Vert_{L^q(0,T;X)}
\le C\left( \int^T_0 t^{-\alpha q(1-\mu)} dt \right)^{\frac{1}{q}}
\Vert (-A)^{\mu}(a_n-a)\Vert\\
\le& C_1\Vert (-A)^{\mu}(a_n-a)\Vert.
\end{align*}
Here we used
$$
\int^T_0 t^{-\alpha q(1-\mu)} dt = \frac{T^{1-\alpha q(1-\mu)}}
{1-\alpha q(1-\mu)} < \infty
$$
by (2.8).  Therefore, $AGa_n \rrrr AGa$ in $L^q(0,T;X)$ as 
$n\to \infty$.

Since $\pppa (Ga_n-a_n) = SAGa_n$ by (2.10), we see that 
$$
\pppa (Ga_n-a_n) \rrrr AGa, \quad
Ga_n - a_n \in \DDD(\pppa).
                             \eqno{(2.19)}
$$
Moreover, 
$$
\Vert (Ga_n - a_n) - (Ga-a)\Vert_{L^q(0,T;X)}
\le \Vert G(a_n - a)\Vert_{L^q(0,T;X)}
+ \Vert a_n-a\Vert_{L^q(0,T;X)} \, \rrrr \, 0
                                          \eqno{(2.20)}
$$
as $n\to \infty$ by means of Proposition 2.1 (ii).
Since $\pppa : L^q(0,T;X) \rrrr L^q(0,T;X)$ is a closed 
operator, the limits (2.19) and (2.20) yield 
$$
Ga - a \in \DDD(\pppa), \quad \mbox{and}\quad
\pppa (G(t)a - a) = AG(t)a.
$$
By means of Proposition 2.1 (ii), we can readily derive that 
$Ga - a \in L^q(0,T;X)$ for all $a\in X$.
Thus the proof of Theorem 2.1 is completed.
$\blacksquare$

\section{Construction of solution operator in the case of 
$a=0$ with $F\ne 0$}

In this section, we consider the case where the initial value is zero and 
the non-homogeneous term is not zero:
$$
\pppa u = Au + F(t),   \quad t>0.              \eqno{(3.1)}
$$
\\
{\bf 3.1. Estimation of $\ABAB \frac{dG}{dt}(t)a$.}
\\
We recall 
$$
G(t)a = \frac{1}{2\pi i}\int_{\Gamma} e^{\la t}\la^{\alpha-1}
(\la^{\alpha} - A)^{-1}a d\la,
$$
where the path $\Gamma \subset \C$ is defined by (2.3).
We are limited to $t>0$, not $z\in \Sigma_{\gamma-\frac{\pi}{2}}$,
although the estimation is possible.

For any fixed $t>0$, we see that 
$$
\frac{\ppp}{\ppp t}(e^{\la t}\la^{\alpha-1}(-A)^{\beta}
(\la^{\alpha}-A)^{-1}a) \in L^{\infty}(\Gamma) \cap L^1(\Gamma)
\quad \mbox{as a function in $\la$.}
$$
In terms of the Lebesgue convergence theorem, 
we can justify the exchange of $\frac{\ppp}{\ppp t}$ and 
$\frac{1}{2\pi i} \int_{\Gamma} \cdots d\la$, so that  
$$
\frac{\ppp}{\ppp t}(G(t)a)
= \frac{1}{2\pi i}\int_{\Gamma} \frac{\ppp}{\ppp t}
(e^{\la t}\la^{\alpha-1}(\la^{\alpha} - A)^{-1}a) d\la.
                                                        \eqno{(3.2)}
$$
Let 
$$
0 \le \beta \le 1.
$$
Then, 
\begin{align*}
& (-A)^{\beta}G'(t)a
= \frac{1}{2\pi i}\int_{\Gamma} e^{\la t}\la^{\alpha}(-A)^{\beta}
(\la^{\alpha} - A)^{-1}a d\la\\
=& \frac{1}{2\pi i}\left( \int_{\Gamma_1} + \int_{\Gamma_2}
+ \int_{\Gamma_3}\right)  e^{\la t}\la^{\alpha}(-A)^{\beta}
(\la^{\alpha} - A)^{-1}a d\la
=: K_1(t) + K_2(t) + K_3(t).
\end{align*}
\\
{\bf Estimation of $K_1(t)$}
We set $\la = \rho e^{-i\gamma}$ with $\rho > \frac{1}{t}$.
In view of Lemma 2.1, using (2.7), we have
\begin{align*}
& \Vert K_1(t)\Vert = \frac{1}{2\pi}\int_{\Gamma_1}
 e^{-\rho t \delta} \vert \la\vert^{\alpha}
\Vert (-A)^{\beta}(\la^{\alpha} - A)^{-1}a\Vert \vert d\la\vert\\
\le& C\int^{\infty}_{\frac{1}{t}} e^{-\rho t\delta} \rho^{\alpha}
\rho^{\alpha(\beta-1)} d\rho \Vert a\Vert
\le C\int^{\infty}_{\frac{1}{t}} e^{-\delta \rho t} \rho^{\alpha\beta} 
d\rho \Vert a\Vert.
\end{align*}
Changing the variables $\rho \mapsto \eta$ in the integral 
$\rho = \frac{\eta}{t}$,
we calculate
$$
  \int^{\infty}_{\frac{1}{t}} e^{-\delta \rho t} \rho^{\alpha\beta} 
d\rho
= t^{-\alpha\beta-1}\int^{\infty}_1 e^{-\delta\eta} 
\eta^{\alpha\beta} d\eta
=: Ct^{-\alpha\beta - 1},
$$
which implies 
$$
\Vert K_1(t)\Vert \le C\frac{1}{t^{\alpha\beta+1}}\Vert a\Vert.
$$
\\
{\bf Estimation of $K_2(t)$}
\\
We set $\la = \frac{1}{t}e^{i\theta}$, where $\theta$ varies from 
$-\gamma$ to $\gamma$.  Then, $d\la = \frac{1}{t}ie^{i\theta}d\theta$ and 
$\vert \la\vert = \frac{1}{t}$, $\vert d\la\vert = \frac{1}{t}d\theta$.
Hence,
\begin{align*}
& \Vert K_2(t)\Vert \le C\int^{\gamma}_{-\gamma}
e^{\cos \theta} \left( \frac{1}{t}\right)^{\alpha}
\vert \la\vert^{\alpha(\beta-1)}\Vert a\Vert \frac{1}{t} d\theta\\
\le & C\int^{\gamma}_{-\gamma}
\frac{1}{t^{\alpha\beta+1}} \Vert a\Vert d\theta
\le Ct^{-\alpha\beta-1} \Vert a\Vert.
\end{align*}

Similarly to $K_1(t)$, we can estimate $K_3(t)$, and we obtain
\\
{\bf Proposition 3.1}
\\
{\it 
Let $0\le \beta \le 1$.  Then
$$
\Vert (-A)^{\beta}\frac{d}{dt}G(t)a\Vert \le Ct^{-\alpha\beta-1}
\Vert a\Vert \quad \mbox{for all $t>0$ and $a \in X$.}
$$
}
\\
\vspace{0.2cm}
\\
{\bf 3.2. Estimation of $\ABAB \frac{d}{dt}(J^{\tau}G(t))$.}
\\
For $0 < \tau \le 1$, we note  
$$
\frac{d}{dt}(J^{\tau}v)(t) 
= \frac{1}{\Gamma(\tau)}\frac{d}{dt}
\int^t_0 (t-s)^{\tau-1}v(s) ds \quad \mbox{for $v\in L^1(0,T;X)$}.
$$

{\bf Remark.}
We remark that $\frac{d}{dt}(J^{\tau}v)(t) = D_t^{1-\tau}v(t)$ is 
the Riemann-Liouville fractional derivative.

Henceforth, for $\alpha_1, \alpha_2 > 0$, we define the Mittag-Leffler 
function by
$$
E_{\alpha_1,\alpha_2}(z):= \sum_{k=0}^\infty \frac{z^k}
{\Gamma(\alpha_1 k + \alpha_2)}
$$
and we know that the radius of convergence of the power series is $\infty$
and $E_{\alpha_1,\alpha_2}(z)$ is an entire function in $z \in \C$
(e.g., Podlubny \cite{Po}).

First we prove
\\
{\bf Lemma 3.1.}
\\
{\it 
Let $0 < \tau \le 1$ and $\la \in \C$.  Then,
$$
\frac{d}{dt}(J^{\tau}e^{\la t})(t)
= t^{\tau-1}E_{1,\tau}(\la t), \quad t>0.
$$
}
\\
{\bf Proof of Lemma 3.1.}
\\
We have
$$
 (J^{\tau}e^{\la t})(t)
= \frac{1}{\Gamma(\tau)}\int^t_0 (t-s)^{\tau-1}
\sumk \frac{\la^ks^k}{k!} ds
= \frac{1}{\Gamma(\tau)}\sumk \frac{\la^k}{k!}
\int^t_0 (t-s)^{\tau-1}s^k ds
$$
$$
= \frac{1}{\Gamma(\tau)}\sumk \frac{\la^k}{k!}
\frac{\Gamma(\tau)\Gamma(k+1)}{\Gamma(\tau+k+1)}
t^{\tau+k}
= \sumk \frac{1}{\Gamma(\tau+k+1)} \la^k t^{k+\tau}
= t^{\tau}E_{1,\tau+1}(\la t),     \eqno{(3.3)}
$$
so that 
\begin{align*}
& \frac{d}{ds}(J^{\tau}e^{\la t})(t)
= \sumk \frac{\la^k(k+\tau)t^{k+\tau-1}}
{\Gamma(\tau+k+1)}
= \sumk \frac{(\tau+k)\la^kt^k}{(\tau+k)\Gamma(\tau+k)}
t^{\tau-1}\\
= & t^{\tau-1}\sum_{k=0}^{\infty} \frac{(\la t)^k}{\Gamma(\tau+k)}
= t^{\tau-1}E_{1,\tau}(\la t).
\end{align*}
Thus the proof of Lemma 3.1 is complete.
$\blacksquare$

Moreover we can show
\\
{\bf Lemma 3.2.}
\\
{\it
We choose $\sigma \in \left( \frac{\pi}{2},\, \pi\right)$ and 
$\tau \in (0, \,1]$ arbitrarily.  Then
$$
\vert E_{1,\tau}(\la t)\vert \le \frac{C}{1+\vert \la t\vert}
\quad \mbox{for all $t\ge 0$ and $\sigma \le \vert \mbox{arg}\, \la\vert 
\le \pi$}.
$$
}

The proof is found as Theorem 1.6 (p.35) in Podlubny \cite{Po}.

We recall that $\frac{\pi}{2} < \gamma < \pi$, so that we can 
choose $\gamma_0$ such that $\frac{\pi}{2} < \gamma_0 < \gamma < \pi$.
Therefore, $\Gamma \subset \{ z\in \C;\, 
\gamma_0 \le \vert \mbox{arg}\, z\vert \le \pi\}$ and
$$
\vert E_{1,\tau}(\la t) \vert \le \frac{C}{1+\vert \la t\vert}
\quad \mbox{for $t\ge 0$ and $\la \in \Gamma$.}
                                                    \eqno{(3.4)}
$$

Now we show
\\
{\bf Lemma 3.3}
\\
{\it
Let $0 \le \beta \le 1$.  Then there exists a constant 
$C=C(\beta) > 0$ such that 
$$
\Vert J^{\beta}G(t)a\Vert \le Ct^{\beta}\Vert a\Vert
$$
for all $t>0$ and $a\in X$.  In particular, 
$J^{\beta}Ga \in L^{\infty}(0,T;X)$.
}
\\
{\bf Proof of Lemma 3.3.}
\\
Similarly to Lemma 2.2, we see
$$
J^{\beta}G(t)a = \frac{1}{2\pi i}\int_{\Gamma}
J_t^{\beta}(e^{\la t})\la^{\alpha-1}(\la^{\alpha}-A)^{-1}a d\la.
$$
Here $J_t^{\alpha}$ means the fractional integral operator 
acting on a function in $t$.
By (3.3), we have
$$
J^{\beta}G(t)a = \frac{t^{\beta}}{2\pi i}\int_{\Gamma}
E_{1,\beta+1}(\la t)\la^{\alpha-1}(\la^{\alpha}-A)^{-1}a d\la.
$$
In terms of (3.4), we obtain
\begin{align*}
& \Vert J^{\beta}G(t)a\Vert 
\le Ct^{\beta} \int_{\Gamma}
\frac{1}{1 + \vert \la t\vert} \vert \la\vert^{\alpha-1}
\frac{1}{\vert \la\vert^{\alpha}} \vert d\la\vert \Vert a\Vert\\
=& Ct^{\beta} \left( \int_{\Gamma_1}
+ \int_{\Gamma_2} + \int_{\Gamma_3}\right)
\frac{1}{1 + \vert \la t\vert} \frac{1}{\vert \la\vert} 
\vert d\la\vert \Vert a\Vert.
\end{align*}
Then,
$$
 \int_{\Gamma_1} \frac{1}{1 + \vert \la t\vert} \frac{1}{\vert \la\vert} 
\vert d\la\vert
= \int^{\infty}_{\frac{1}{t}} \frac{1}{1+\rho t}\frac{1}{\rho} d\rho
= \int^{\infty}_1 \frac{1}{1+\eta} \frac{1}{\eta} d\eta
=: C_0.
$$
Here the constant $C_0>0$ is independent of $t>0$, and we used 
the change of variables $\eta:= \rho t$: $\rho \mapsto \eta$ for the last 
integration.

Similarly we can estimate $\int_{\Gamma_3} \frac{1}{1 + \vert \la t\vert} 
\frac{1}{\vert \la\vert} \vert d\la\vert$.
Finally 
$$
\int_{\Gamma_2} \frac{1}{1 + \vert \la t\vert} \frac{1}{\vert \la\vert} 
\vert d\la\vert
= \int^{\gamma}_{-\gamma} \frac{1}{2} t \frac{1}{t} d\theta = \gamma,
$$
where $\la\in \Gamma_2$ is parametrized by $\la = \frac{1}{t}e^{i\theta}$
with $\theta \in (-\gamma,\gamma)$.

Therefore, we can choose a constant $C>0$ such that 
$\Vert J^{\beta}G(t)a\Vert \le Ct^{\beta}\Vert a\Vert$ for all 
$t>0$ and $a\in X$.  Thus the proof of Lemma 3.3 is complete.
$\blacksquare$
\\

Now we start to estimate $(-A)^{\beta}\frac{d}{dt}J^{\tau}G$:
$$
 (-A)^{\beta}\frac{d}{dt}J^{\tau}G(t)a 
= \frac{1}{2\pi i}\int_{\Gamma} \frac{d}{dt}J^{\tau}(e^{\la t})\la^{\alpha-1}
(-A)^{\beta}(\la^{\alpha} - A)^{-1}a d\la.
                                                    \eqno{(3.5)}
$$
In this equality, the verification of the exchange of $\frac{d}{dt}J^{\tau}$ 
and $\int_{\Gamma} \cdots d\la$, is similar to (3.2), and we omit the
details.  Therefore, in terms of Lemma 3.1 and (3.3), we have 
\begin{align*}
& \Vert (-A)^{\beta} \frac{d}{dt}J^{\tau}G(t)a\Vert
\le \frac{1}{2\pi}\int_{\Gamma} t^{\tau-1} \vert E_{1,\tau}(\la t)\vert
\vert \la\vert^{\alpha-1}\Vert (-A)^{\beta}(\la^{\alpha} - A)^{-1}a\Vert
\vert d\la\vert                               \\
\le &Ct^{\tau-1} \int_{\Gamma}\frac{1}{1+\vert \la t\vert}
\vert \la\vert^{\alpha-1}\vert \la\vert^{\alpha(\beta-1)} \vert d\la\vert
\Vert a\Vert
= Ct^{\tau-1}\int_{\Gamma} \frac{\vert \la\vert^{\alpha\beta-1}}
{1+\vert \la t\vert} \vert d\la\vert \Vert a\Vert\\
= & Ct^{\tau-1}\Vert a\Vert \left(\int_{\Gamma_1}
+ \int_{\Gamma_2} + \int_{\Gamma_3}\right) 
\frac{\vert \la\vert^{\alpha\beta-1}}
{1+\vert \la t\vert} \vert d\la\vert
=: J_1(t) + J_2(t) + J_3(t).
\end{align*}

{\bf Estimate of $J_1(t)$}.
We set $\la = \rho e^{-i\gamma}$ with $\rho > \frac{1}{t}$.
In view of Lemma 2.1, also by the change $\eta:= \rho t$ and 
$\alpha\beta - 1 < 0$, we have
\begin{align*}
& \int_{\Gamma_1} \frac{\vert \la\vert^{\alpha\beta-1}}
{1+ \vert \la t\vert} \vert d\la \vert
= \int^{\infty}_{\frac{1}{t}} \frac{\rho^{\alpha\beta-1}}
{1 + \rho t} d\rho
= \int^{\infty}_1 \left( \frac{\eta}{t}\right)^{\alpha\beta-1}
\frac{1}{1+\eta}\frac{1}{t} d\eta\\
=& t^{-\alpha\beta}\int^{\infty}_1 \frac{1}{\eta^{1-\alpha\beta}}
\frac{1}{1+\eta} d\eta
\le t^{-\alpha\beta}\int^{\infty}_1 \eta^{\alpha\beta-2} d\eta
=: Ct^{-\alpha\beta}.
\end{align*}
Therefore,
$$
\Vert J_1(t)\Vert \le Ct^{\tau-1}t^{-\alpha\beta} \Vert a\Vert
= Ct^{(\tau-\alpha\beta)-1}\Vert a\Vert.
$$

{\bf Estimation of $J_2(t)$.}
We set $\la = \frac{1}{t}e^{i\theta}$, where $\theta$ varies from 
$-\gamma$ to $\gamma$.  Then, $d\la = \frac{1}{t}ie^{i\theta}d\theta$ and 
$\vert \la\vert = \frac{1}{t}$, $\vert d\la\vert = \frac{1}{t}d\theta$.
Consequently, 
$$
\int_{\Gamma_2} \frac{\vert \la\vert^{\alpha\beta-1}}
{1+\vert \la t\vert}\vert d\la\vert
= \int_{-\gamma}^{\gamma} \left( \frac{1}{t}\right)^{\alpha\beta-1}\frac{1}{2}
\frac{1}{t} d\theta
= \frac{1}{2}t^{-\alpha\beta} (2\gamma).
$$
Hence, 
$$
\Vert J_2(t)\Vert \le Ct^{\tau-1}t^{-\alpha\beta} \Vert a\Vert
= Ct^{(\tau-\alpha\beta)-1}\Vert a\Vert.
$$
As for $J_3(t)$, the estimation is essentially same as $J_1(t)$.
Thus we have proved
\\
{\bf Proposition 3.2.}
\\
{\it
Let $0 \le \beta \le 1$ and $0 < \tau \le 1$.  Then, we have
$$
\frac{d}{dt}J^{\tau}Ga\in L^1(0,T;X), \quad  
\left\Vert (-A)^{\beta} \frac{d}{dt}J^{\tau}G(t)a\right\Vert 
\le Ct^{(\tau-\alpha\beta)-1}\Vert a\Vert, \quad t>0.
$$
In particular, choosing $\tau = \alpha$, we have
$$
\left\Vert (-A)^{\beta} \frac{d}{dt}J^{\alpha}G(t)a\right\Vert 
\le Ct^{\alpha(1-\beta)-1}\Vert a\Vert \quad \mbox{for all 
$a \in X$ and $t>0$}.
$$
}
\\

We set 
$$
K(t)a: = \frac{d}{dt}J^{\alpha}G(t) a, \quad a\in X.   \eqno{(3.6)}
$$

The following lemma can be proved by similar estimation in this section,
and is not used in 
Sections 3 - 5 but is useful in Section 6.
\\
{\bf Lemma 3.4.}
{\it
(i)
$$
\Vert AK(t)a - AK(s)a\Vert \le Ct^{-\alpha}(t^{\alpha-1}
- s^{\alpha-1})\Vert a\Vert
$$
for all $0<t \le s<T$ and $a \in X$.
\\
(ii)
$$
\Vert AJ^{\alpha}G(t)a \Vert \le C\Vert a\Vert
$$
for all $0\le t \le T$ and $a \in X$.
}
\\
{\bf Proof of Lemma 3.4.}
\\
Fixing $\delta > 0$ sufficiently small,
we define a path $\Gamma_{\delta}:= 
\Gamma_\delta^1 \cup \Gamma_\delta^2 \cup \Gamma_\delta^3$ in $\C$ by 
$$
\Gamma_\delta^1 := \{ \rho e^{-i\gamma};\, \rho > \delta\}, \quad
\Gamma_\delta^2 := \{ \delta e^{i\theta};\, -\gamma \le \theta \le 
\gamma\}, \quad \Gamma_\delta^3 := \{ \rho e^{i\gamma};\, \rho > \delta\}, 
$$
and $\Gamma_{\delta}$ is directed from $\infty e^{-i\gamma}$ to 
$\infty e^{i\gamma}$.  By Condition ($\AAA$), we note that 
the integral in (2.2) is invariant with repect to $\Gamma_{\delta}$
with small $\delta>0$.

By means of (3.5) and Lemma 3.1, we have
\begin{align*}
& I:= AK(t)a - AK(s) a 
=  A\frac{d}{dt}J^{\alpha}G(t)a  - A\frac{d}{ds}J^{\alpha}G(s)a\\
=& \frac{1}{2\pi i} \int_{\Gamma_{\delta}}
\left( \frac{d}{dt}J^{\alpha}e^{\la t}
- \frac{d}{ds}J^{\alpha}e^{\la s}\right) \la^{\alpha-1}A
(\la^{\alpha} - A)^{-1}a d\la \\
=& \frac{1}{2\pi i} \int_{\Gamma_{\delta}}
 ( t^{\alpha-1}E_{1,\alpha}(\la t) - s^{\alpha-1}E_{1,\alpha}(\la s))
 \la^{\alpha-1}A(\la^{\alpha} - A)^{-1}a d\la.
\end{align*}

Since 
$$
\frac{d}{dt}(t^{\alpha-1}E_{1,\alpha}(\la t)) = t^{\alpha-2}E_{1,\alpha-1}(\la t)
$$
(e.g., (1.83) (p.22) in \cite{Po}), we have
$$
 t^{\alpha-1}E_{1,\alpha}(\la t) - s^{\alpha-1}E_{1,\alpha}(\la s)
= \int^t_s \xi^{\alpha-2}E_{1,\alpha-1}(\la \xi) d\xi,
$$
and so 
\begin{align*}
& I = \frac{1}{2\pi i} \int_{\Gamma_{\delta}}
 \left( \int^t_s \xi^{\alpha-2}E_{1,\alpha}(\la \xi)) d\xi\right)
 \la^{\alpha-1}A(\la^{\alpha} - A)^{-1}a d\la\\
=& \frac{1}{2\pi i} \left(\int_{\Gamma_{\delta}^1}
 + \int_{\Gamma_{\delta}^2} + \int_{\Gamma_{\delta}^3}\right)
 \left( \int^t_s \xi^{\alpha-2}E_{1,\alpha}(\la \xi)) d\xi\right)
 \la^{\alpha-1}A(\la^{\alpha} - A)^{-1}a d\la 
=: I_1 + I_2 + I_3.
\end{align*}
Therefore, using Condition ($\AAA$), we obtain
$$
\Vert I_k\Vert \le C\left( \int_{\Gamma_{\delta}^k} 
 \left( \int^s_t \xi^{\alpha-2}\vert E_{1,\alpha}(\la \xi)\vert d\xi\right)
 \vert \la\vert^{\alpha-1}\vert d\la\vert \right) \Vert a\Vert, \quad k=1,2,3.
$$
On the other hand, Theorem 1.6 (p.35) in \cite{Po} yields
$$
\vert E_{1,\alpha-1}(z)\vert \le \frac{C}{1 + \vert z\vert}
\quad \mbox{for all $z \in \C$ satisfying
$\frac{1}{2}\pi < \vert \mbox{arg}\, z\vert \le \pi$}.
$$
Hence, since $\frac{1}{2}\pi < \vert \mbox{arg}\, (\la\xi)\vert \le \pi$
for $\xi \ge 0$ and $\la \in \Gamma_{\delta}^1 \cup \Gamma_{\delta}^3$,
we see
$$
\vert E_{1,\alpha-1}(\la \xi)\vert \le \frac{C}{1 + \vert \la\xi\vert}
\quad \mbox{for all $\la \in \Gamma_{\delta}^1 \cup 
\Gamma_{\delta}^3$ and $\xi \in [0,T]$}.
$$
On the other hand, if $\la \in \Gamma_{\delta}^2$, then 
$\vert \la \xi\vert = \delta \xi \le \delta T$ and so
$$
\vert E_{1,\alpha-1}(\la\xi)\vert \le C_1 = C_1(\delta,T)
\quad \mbox{for all $\la \in \Gamma_{\delta}^2$ and $\xi \in [0,T]$}.
$$
Therefore, we can find a sufficiently large constant 
$C_2 = C_2(\delta,T)>0$ such that 
$$
\vert E_{1,\alpha-1}(\la\xi)\vert \le \frac{C_2}{1 + \vert \la\vert \xi}
\le \frac{C_2}{1 + \vert \la\vert t}
\quad \mbox{for $\la \in \Gamma_{\delta}, \, \xi \in [t,s]$}.   \eqno{(3.7)}
$$
We estimate $I_1, I_2, I_3$ separately.
\\
{\bf Estimation of $I_1$.}
We have
\begin{align*}
& \Vert I_1\Vert 
\le C_3\int_{\Gamma_{\delta}^1} \left( \int^s_t \xi^{\alpha-2} d\xi\right)
\frac{\vert \la\vert^{\alpha-1}}{1 + \vert \la \vert t}\vert d\la\vert 
\Vert a\Vert\\
= & C_3\Vert a\Vert \int_{\Gamma_{\delta}^1} \frac{1}{1-\alpha}
(t^{\alpha-1} - s^{\alpha-1}) \frac{\vert \la\vert^{\alpha-1}}
{1 + \vert \la \vert t}\vert d\la\vert
= C_4\Vert a\Vert (t^{\alpha-1} - s^{\alpha-1}) 
\int^{\infty}_{\delta} \frac{\rho^{\alpha-1}}{1 + \rho t} d\rho.
\end{align*}
In the last equality, we used the change of the variables 
$\la = \rho e^{-i\gamma}$ with $\rho > \delta$.
Again the change of the variables $\rho \mapsto \eta$ by 
$\eta = t\rho$, we have 
$$
\int^{\infty}_{\delta} \frac{\rho^{\alpha-1}}{1 + \rho t} d\rho
= t^{-\alpha}\int^{\infty}_{\delta t} \frac{\eta^{\alpha-1}}{1 + \eta} 
d\eta
\le t^{-\alpha}\int^{\infty}_0 \frac{1}{\eta^{1-\alpha}(1 + \eta)}d\eta
=: C_5t^{-\alpha}.
$$
Hence, $\Vert I_1\Vert \le C_5t^{-\alpha}(t^{\alpha-1} - s^{\alpha-1})
\Vert a\Vert$ for $0<t\le s < T$,
where $C_5 = C_5(\delta,T)>0$.
In the same way, we can obtain the same estimate $I_3$.
\\
{\bf Estimation of $\Vert I_2\Vert$.}
\\
For $\la\in \Gamma^2_{\delta}$, we can give
$\la = \delta e^{i\theta}$ where $\theta: -\gamma \mapsto \gamma$, and so
$\vert \la\vert = \delta$ and $\vert d\la\vert = \delta d\theta$.
Therefore,
\begin{align*}
& \Vert I_2\Vert \le C\left( \int^{\gamma}_{-\gamma}
\left( \int^s_t \xi^{\alpha-2}\frac{C_2}{1+\delta t} d\xi \right)
\times \delta^{\alpha-1}\delta d\theta \right)\Vert a\Vert\\
\le& CC_2 \delta^{\alpha}\Vert a\Vert \int^{\gamma}_{-\gamma}
\left( \int^s_t \xi^{\alpha-2} d\xi\right) d\theta \\
\le & C_6(t^{\alpha-1} - s^{\alpha-1}) \Vert a\Vert 
\le C_7 t^{-\alpha} (t^{\alpha-1} - s^{\alpha-1}),
\end{align*}
where we used $T^{-\alpha} \le t^{-\alpha}$ for $0 <t \le T$.
Thus the proof of (i) is complete.
\\
{\bf Proof of (ii).}
\\
By (4.2) and (3.3), we obtain
\begin{align*}
& AJ^{\alpha}G(t)a = \frac{1}{2\pi i}
\int_{\Gamma_\delta} J^{\alpha}(e^{\la t}) \la^{\alpha-1}
A (\la^{\alpha} - A)^{-1} a d\la \\
=& \frac{1}{2\pi i}
\int_{\Gamma_\delta} t^{\alpha} E_{1,\alpha+1}(\la t)
 \la^{\alpha-1}A(\la^{\alpha}-A)^{-1} a d\la.
\end{align*}
Therefore,
\begin{align*}
& \Vert AJ^{\alpha}G(t)a \Vert 
\le C\int_{\Gamma_\delta} t^{\alpha} \vert E_{1,\alpha+1}(\la t)\vert
\vert \la\vert^{\alpha-1} \vert d\la\vert \Vert a\Vert\\
=& C\left( \int_{\Gamma_\delta^1} + \int_{\Gamma_\delta^2}
    + \int_{\Gamma_\delta^3}\right) 
 t^{\alpha} \vert E_{1,\alpha+1}(\la t)\vert
\vert \la\vert^{\alpha-1} \vert d\la\vert \Vert a\Vert
=: S^1 + S^2 + S^3.
\end{align*}
Using (3.7) and {\it Condition ($\AAA$)}, we see
\begin{align*}
& \Vert S^1\Vert \le C\int_{\Gamma_\delta^1}
t^{\alpha}\frac{1}{1 + \vert \la \vert t} \vert \la\vert^{\alpha-1}
\vert d\la \vert \Vert a\Vert 
= C\Vert a\Vert\int_{\delta}^{\infty} t^{\alpha}\frac{1}{1 + \rho t} 
\rho^{\alpha-1}  d\rho  \\
= & C\Vert a\Vert \int_{\delta t}^{\infty} \eta^{\alpha-1}\frac{1}{1 + \eta} 
d\eta \le C\left( \int^{\infty}_0 \frac{d\eta}{\eta^{1-\alpha})1+\eta)}
\right) \Vert a\Vert \le C_7\Vert a\Vert.
\end{align*}
Similarly we can see $\Vert S^3\Vert \le C_7\Vert a\Vert$ and 
make estimation of $\Vert S^2\Vert$.
Thus the proof of Lemma 3.4 is complete.
$\blacksquare$
\\
\vspace{0.2cm}
\\ 
{\bf 3.3. Equation to be satisfied.}
\\
We show
\\
{\bf Lemma 3.5.}
\\
{\it
Let $F \in C^{\infty}_0(0,T;\DDD(A))$.  Then
$$
\int^t_0 \left( \frac{d}{dt}J^{\alpha}G\right)(t-s)F(s) ds
\in W_{\alpha,q}(0,T;X)
$$
and
$$
\pppa \int^t_0 \left( \frac{d}{dt}J^{\alpha}G\right)(t-s)F(s) ds
= \int^t_0 G(t-s)F'(s) ds, \quad t>0.                
$$
}
\\
{\bf Proof of Lemma 3.5.}
\\
By Lemma 3.3, we see that $J^{\alpha}Ga \in L^{\infty}(0,\infty;X)$ for 
all $a \in X$.  Hence, for sufficiently small $\delta>0$, by integration 
by parts, we have
\begin{align*}
& \int^{t-\delta}_0 \left( \frac{d}{dt}J^{\alpha}G\right)(t-s)F(s) ds
= -\int^{t-\delta}_0 \frac{d}{ds}(J^{\alpha}G(t-s))F(s) ds \\
= -& [J^{\alpha}G(t-s)F(s)]^{s=t-\delta}_{s=0}
+ \int^{t-\delta}_0 J^{\alpha}G(t-s)F'(s) ds.
\end{align*}
Since $F \in C^{\infty}_0(0,T;\DDD(A))$, we have $F(0) = 0$. 
Lemma 3.3 yields
$$
\Vert J^{\alpha}G(\delta)F(t-\delta)\Vert \le C\delta^{\alpha}
\Vert F(t-\delta)\Vert 
\le C\Vert F\Vert_{C([0,T];X)}\delta^{\alpha} \rrrr 0
$$
as $\delta \to 0$.  Consequently,
$$
\int^t_0 \left( \frac{d}{dt}J^{\alpha}G\right)(t-s)F(s) ds
= \int^t_0 J^{\alpha}G(t-s)F'(s) ds, \quad 0<t<T.       \eqno{(3.8)}
$$

Next we will calculate:
$$
J^{\alpha}\left( \int^t_0 G(t-s)F'(s) ds\right)
= \frac{1}{\Gamma(\alpha)} \int^t_0 (t-s)^{\alpha-1}
\left( \int^s_0 G(s-\xi)F'(\xi) d\xi \right) ds.
$$
Exchanging the orders of the integrals:
$$
\int^t_0 \left( \int^s_0 \cdots d\xi\right) ds
= \int^t_0 \left( \int^t_\xi \cdots ds\right) d\xi,
$$
we obtain
\begin{align*}
& \frac{1}{\Gamma(\alpha)} \int^t_0 (t-s)^{\alpha-1}
\left( \int^s_0 G(s-\xi)F'(\xi) d\xi \right) ds
= \frac{1}{\Gamma(\alpha)} \int^t_0 
\left( \int^t_\xi (t-s)^{\alpha-1}G(s-\xi) ds\right) F'(\xi) d\xi\\
=& \frac{1}{\Gamma(\alpha)} \int^t_0 
\left( \int^{t-\xi}_0 (t-\xi-\eta)^{\alpha-1}G(\eta) d\eta\right) 
F'(\xi) d\xi.
\end{align*}
For the last equality, we changed the variables: $s \mapsto \eta$ by 
$\eta:= s-\xi$.

Since $\frac{1}{\Gamma(\alpha)} \int^{t-\xi}_0 
(t-\xi-\eta)^{\alpha-1}G(\eta) d\eta = (J^{\alpha}G)(t-\xi)$, 
by means of (3.8), we reach 
$$
\int^t_0 \left( \frac{d}{dt}J^{\alpha}G\right)(t-s) F(s) ds
= J^{\alpha} \left(\int^t_0 G(t-s) F'(s) ds\right), \quad 
0<t<T.                                  \eqno{(3.9)}
$$
Moreover Proposition 2.1 yields
$$
\left\Vert G(t-s)F'(s) ds\right\Vert
\le \int^t_0 \Vert G(t-s)\Vert \Vert F'(s)\Vert ds 
\le C\max_{0\le s\le T} \Vert F'(s)\Vert,
$$
so that 
$$
\int^t_0 G(t-s)F'(s) ds \in L^q(0,T;X).
$$
Therefore (3.9) and the definition of $\pppa$ imply
$$
J^{\alpha}\int^t_0 G(t-s) F'(s) ds \in W_{\alpha,q}(0,T;X),
$$
that is,
$$
\int^t_0 \left( \frac{d}{dt}J^{\alpha}G\right)(t-s) F(s) ds
\in W_{\alpha,q}(0,T;X).
$$
Thus the proof of Lemma 3.5 is complete.
$\blacksquare$
\\

We set 
$$
W(F)(t):= \int^t_0 \left( \frac{d}{dt}J^{\alpha}G\right)(t-s)F(s) ds
= - \int^t_0 \frac{d}{ds}(J^{\alpha}G(t-s))F(s) ds
$$
for $0<t<T$ and $F \in C^{\infty}_0(0,T;\DDD(A))$.
Then,
\\
{\bf Lemma 3.6.}
\\
{\it
Let $F \in C^{\infty}_0(0,T;\DDD(A))$.  Then,
$W(F) \in W_{\alpha,q}(0,T;X)$, $AW(F) \in L^q(0,T;X)$
and
$$
\pppa W(F)(t) =  AW(F)(t) + F(t) 
\quad \mbox{for $0<t<T$}.
$$
}
\\
{\bf Proof of Lemma 3.6.}
\\
Choosing $\delta > $ sufficiently small, 
we set
$$
W_{\delta}(F)(t) = -\int^{t-\delta}_0 \left(\frac{d}{ds}
(J^{\alpha}G)(t-s)\right) F(s) ds.
$$
Since $F(t)\in \DDD(A)$, we obtain
$$
A\left( \frac{d}{ds} (J^{\alpha}G)(t-s)\right) F(s)
= \left(\frac{d}{ds}(J^{\alpha}G)(t-s)\right) AF(s)
\in L^1_s(0,T;X).
$$
Noting that $A$ is a closed operator in $X$, by an argument similar to 
Lemma 2.2, we 
see
$$
AW_{\delta}(F)(t) = -\int^{t-\delta}_0 \frac{d}{ds}
(J^{\alpha}G)(t-s) AF(s) ds
= -\int^{t-\delta}_0 \left(\frac{d}{ds}(J^{\alpha}AG)(t-s) \right) F(s) ds.
$$
Theorem 2.1 yields 
$$
AG(t-s)a = \pppa (G(t-s)a-a) \quad \mbox{in $X$ for $0<s<t$ and all 
$a\in \DDD(A)$}
$$
and
$$
J^{\alpha}AG(t-s)F(s) = J^{\alpha}\pppa (G(t-s) -1)F(s).
$$
Since $F(t)\in \DDD(A)$, we can apply (2.9) and 
$(G(t-s)-1)F(s) \in \DDD(\pppa)$ as a function in $t$.
Consequently,
$$
J^{\alpha}\pppa (G(t-s)-1) = G(t-s) -1,
$$
which implies
$$
J^{\alpha}AG(t-s)F(s) = (G(t-s)-1)F(s).
$$
Hence,
$$
AW_{\delta}(F)(t) = -\int^{t-\delta}_0 \frac{d}{ds}(G(t-s) - 1)F(s) ds,
\quad 0<t<T.
$$

Similarly to the proof of (3.8), the integration by parts implies
\begin{align*}
& AW_{\delta}(F)(t) 
= -[(G(t-s)-1)F(s)]^{s=0}_{s=t-\delta}
+ \int^{t-\delta}_0 (G(t-s)-1)F'(s) ds\\
=& (G(\delta)-1)F(t-\delta)
+ \int^{t-\delta}_0 G(t-s)F'(s) ds
- \int^{t-\delta}_0 F'(s) ds.
\end{align*}
On the other hand,
\\
{\bf Lemma 3.7.}
{\it
$$
\lim_{t\to 0} \Vert G(t)a - a\Vert = 0
$$
for all $a \in \DDD(A)$.
}

Lemma 3.7 is a special case of Proposition 4.1 stated in 
Section 4, and for convenience, we prove the lemma at the end of this section.
We continue the proof of Lemma 3.6.  We can show  
$$
\lim_{\delta\to 0} \Vert (G(\delta)-1)F(t-\delta)\Vert 
= 0 \quad \mbox{for each $0 < t< T$.}
$$
Indeed, 
\begin{align*}
& \Vert (G(\delta) - 1)F(t-\delta)\Vert
= \Vert (G(\delta) - 1)F(t) + (G(\delta) - 1)(F(t-\delta) - F(t))\Vert\\
\le& \Vert (G(\delta) - 1)F(t) \Vert 
+ \Vert (G(\delta) - 1)(F(t-\delta) - F(t))\Vert.
\end{align*}
Let $t>0$ be arbitrarily fixed.
Then $\lim_{\delta\to 0}\Vert (G(\delta) - 1)F(t)\Vert = 0$ by 
Lemma 3.7.  Moreover, Proposition 2.1 (ii) yields
\begin{align*}
&\Vert G(\delta)-1\Vert \Vert F(t-\delta) - F(t)\Vert 
\le C(\Vert G(\delta)\Vert + 1)\Vert F(t-\delta) - F(t)\Vert \\
\le& C\Vert F(t-\delta) - F(t)\Vert.
\end{align*}
Since $t>0$, the limit of the right-hand side is zero as 
$\delta \to 0$.
\\
Hence $\lim_{\delta\to 0} \Vert (G(\delta) - 1)F(t-\delta)\Vert = 0$
for any fixed $t>0$.
$\blacksquare$

Therefore,
\begin{align*}
& AW(F)(t) = \lim_{\delta\to 0} AW_{\delta}(F)(t)\\
=& \lim_{\delta\to 0} (G(\delta)-1)F(t-\delta)
+ \lim_{\delta \to 0}\int^{t-\delta}_0 G(t-s)F'(s) ds
- \lim_{\delta \to 0} (F(t-\delta) - F(0)),
\end{align*}
and so
$$
AW(F)(t) = \int^t_0 G(t-s)F'(s) ds - F(t) \quad 
\mbox{for $0 < t < T$}.                     \eqno{(3.10)}
$$
On the other hand, the application of Lemma 3.5 to (3.10)
yields
$$
AW(F)(t) = \pppa W(F)(t) - F(t), \quad 0<t<T.
$$
Thus the proof of Lemma 3.6 is complete.
$\blacksquare$
\\

We conclude this subsection with 
\\
{\bf Proof of Lemma 3.7.}
\\
Theorem 2.1 and $a \in \DDD(A)$ yield
$\pppa (G(t)a - a) = AG(t)a = G(t)Aa$.
Therefore, since (2.9) implies $Ga - a\in \DDD(\pppa)$, we obtain
$$
G(t)a - a = J^{\alpha}\pppa (G(t)a - a) 
= J^{\alpha}G(t)Aa.
$$
Applying Lemma 3.3, we see
$$
\Vert G(t)a - a\Vert \le Ct^{\alpha}\Vert Aa\Vert
\quad \mbox{for $t>0$.}
$$
Hence, $\lim_{t\to 0} \Vert G(t)a - a\Vert = 0$, which completes the 
proof of Lemma 3.7.
$\blacksquare$
\section{Solution formula and the unique existence of strong solution}

{\bf 4.1. Existence of solution.}
\\
In this subsection, we summarize the results obtained in Sections 2 and 3, and 
add other facts which are proved readily to state Theorem 4.1 as the main 
result for this article.
 
We consider 
$$
\pppa (u(t)-a) = Au(t) + F(t), \quad 0<t<T.          \eqno{(4.1)}
$$
{\bf Definition 4.1.}
\\
{\it 
Let $T>0$ be given.
We call $u=u(t)$ a {\bf strong solution} to (4.1) if $u\in L^q(0,T;\DDD(A))$ 
satisfies (4.1) and $u - a\in W_{\alpha,q}(0,T;X)$.
}
\\

We can state
\\
{\bf Theorem 4.1.}
\\
{\it
We define  
$$
G(t)a := \frac{1}{2\pi i}\int_{\Gamma} e^{\la t}\la^{\alpha-1}
(\la^{\alpha}-A)^{-1} a d\la, \quad t>0, \, a\in X,             \eqno{(4.2)}
$$
where $\Gamma \subset \C$ is defined by (2.3).
We set 
$$
K(t)a:= \frac{d}{dt}J^{\alpha}G(t)a, \quad t>0, \, a\in X.
                                                   \eqno{(4.3)}
$$
Let 
$$
1\le q \le \infty.
$$
\\
(1) Let 
$$
1 - \frac{1}{q\alpha} < \mu < 1             \eqno{(4.4)}
$$
and let 
$$
a \in \DDD((-A)^{\mu}) \quad \mbox{and}\quad
F \in L^q(0,T; \DDD((-A)^{\ep}))                 \eqno{(4.5)}
$$
with some $\ep > 0$.
Then there exists a strong solution $u$ to (4.1), that is,
$u\in L^q(0,T;\DDD(A))$ and $u-a \in W_{\alpha,q}(0,T;X)$.
Moreover, we can find a constant $C>0$ such that 
$$
\Vert u\Vert_{L^q(0,T;\DDD(A))} + \Vert u-a \Vert_{W_{\alpha,q}(0,T;X)}
\le C(\Vert a\Vert_{\DDD((-A)^{\mu})} 
+ \Vert F\Vert_{L^q(0,T;\DDD((-A)^{\ep}))})      \eqno{(4.6)}
$$
for all $a\in \DDD((-A)^{\mu})$ and $F \in L^q(0,T;\DDD((-A)^{\ep}))$.
\\
For $a\in \DDD((-A)^{\mu})$ and $F \in L^q(0,T;\DDD((-A)^{\ep}))$, 
the strong solution is given by
$$
u(t) = G(t)a + \int^t_0 K(t-s)F(s) ds, \quad 0<t<T.     \eqno{(4.7)}
$$
\\
(2) For each $a\in X$, the function $G(z)a$ in $z>0$ can be holomorphically 
extended to $z\in \Sigma_{\gamma-\frac{\pi}{2}}$ and
$$
L(Ga)(\la) = \la^{\alpha-1}(\la^{\alpha}-A)^{-1}a \quad 
\mbox{for $\la \in \Sigma_{\gamma}$}.
$$
\\
(3) For $0\le \beta \le 1$, there exists a constant $C=C(\beta)>0$ 
such that 
$$
\Vert \ABAB G(t)a\Vert \le Ct^{-\alpha\beta}\Vert a\Vert \quad 
\mbox{for all $t>0$ and $a \in X$}
                                                \eqno{(4.8)}
$$
and
$$
\Vert \ABAB K(t)a\Vert \le Ct^{\alpha(1-\beta)-1}\Vert a\Vert \quad 
\mbox{for all $t>0$ and $a \in X$.}
                                                \eqno{(4.9)}
$$
}

Theorem 4.1 asserts that for the existence of solution 
$u \in L^q(0,T;\DDD(A))$ satisfying 
$u-a \in W_{\alpha,q}(0,T;X)$, we have to assume a stronger 
regularity $F \in L^q(0,T;\DDD((-A)^{\ep}))$ with some $\ep > 0$.
For $q\ne 2$ and a general Banach space, 
we do not know whether we can choose 
$\ep=0$, that is, $F \in L^q(0,T;X)$ is sufficient.
In particular, in the case of $a=0$, we call the maximum regularity 
if $F \in L^q(0,T;X)$ implies that 
$u\in L^q(0,T;\DDD(A)) \cap W_{\alpha,q}(0,T;X)$. 
The maximum regularity is known for the case of $q=2$ and a Hilbert space 
$X$ (e.g., Kubica, Ryszewska and Yamamoto \cite{KRY},
Sakamoto and Yamamoto \cite{SY}, Yamamoto \cite{Y22}, Zacher \cite{Za}). 

In Theorem 4.1 (2), we emphasize that the constant $C>0$ is uniform for 
all $t>0$.  Thus we can state
\\
{\bf Corollary 4.1.}
\\
{\it
We assume {\it Condition ($\AAA$)} and $0 \in \rho(A)$.
For $a\in \DDD((-A)^{\mu})$, let $u(t)$ be the solution to 
$$
\pppa (u-a) = Au(t), \quad t>0.
$$
Then we have a decay estimate 
$$
\Vert u(t)\Vert 
\le Ct^{-\alpha}\Vert a\Vert \quad \mbox{for all $t>0$}.
$$
}

As for a special case $X = L^2(\OOO)$, see
Vergara and Zacher \cite{VZ}, Kubica, Ryszewska and Yamamoto \cite{KRY},
for example.  
Vergara and Zacher \cite{VZ} discusses the case of symmetric $A$, 
including a time-fractional $p$-Laplacian equation and a
porus medium equation, while \cite{KRY} is concerned with
linear but not necessarily
symmetric $A$.  Both works are concerned with the cases where 
the coefficients depend on $x$ and $t$.

Our decay estimate is proved directly as by-product 
within a general setting including $X=L^p(\OOO)$.  
\\

Theorem 4.1 does not conclude the continuity in $t$ of the solution, so that
the initial condition in the sense $\lim_{t\to\infty} \Vert u(t) - a\Vert = 0$
is not clear.  In the case of $F=0$, we can prove the following proposition.
As for more regularity, we will study in Subsection 6.2.
\\
{\bf Proposition 4.1.}
$$
\lim_{t'\to t} \Vert G(t')a - G(t)a \Vert = 0
\quad \mbox{for all $0\le t \le T$ and $a\in X$}.
$$

Here we interpret $\lim_{t' \to t} = \lim_{t'\downarrow t}$ for 
$t=0$ and $\lim_{t' \to t} = \lim_{t'\uparrow t}$ for 
$t=T$.
\\
\vspace{0.2cm}
\\
{\bf Proof of Theorem 4.1.}
\\
Part (2) directly follows from Theorem 2.6.1 (p.84) in \cite{ABHN}.
The proof of part (3) of the theorem is already completed by 
Propositions 2.1 and 3.1.
The case $F\equiv 0$ in part (1) is also complete.
Thus the rest is the proof of part (1) in the case of $a=0$ and $F \in L^q(0,T;
\DDD((-A)^{\ep}))$.

We note that (1) with $a=0$ and $F\in C^{\infty}_0(0,T;\DDD(A))$
was also proved by Lemma 3.4:
\begin{align*}
& W(F)(t):= \int^t_0 K(t-s)F(s) ds \in W_{\alpha,q}(0,T;X) \cap
L^q(0,T;X),\\
& \pppa W(G)(t) = AW(F)(t) + F(t), \quad 0<t<T \quad
\mbox{for $F \in C^{\infty}_0(0,T;\DDD(A))$}.
\end{align*}
Now let $F \in L^q(0,T;\DDD(\AEAE))$ be given arbitrarily.
Then, by the mollifier in $t$ (e.g., Adams \cite{Ad}) and the 
density of $\DDD(A)$ in $\DDD(\AEAE)$, we can find a sequence 
$F_n \in C^{\infty}_0(0,T;\DDD(A))$, $n\in \N$ such that 
$\Vert F_n-F\Vert_{L^q(0,T;\DDD(\AEAE))} \rrrr 0$ as 
$n\to \infty$.
In terms of (4.9), we have
\begin{align*}
& \Vert AW(F_n)(t) - AW(F)(t)\Vert_X
= \left\Vert \int^t_0 (-A)^{1-\ep}(-A)^{-\ep} K(t-s)(F-F_n)(s) ds
\right\Vert_X\\
=& \left\Vert \int^t_0 (-A)^{1-\ep}K(t-s)(-A)^{-\ep} (F-F_n)(s) ds
\right\Vert_X
\le C\int^t_0 (t-s)^{\alpha\ep-1} \Vert \AEAE (F-F_n)(s)\Vert_X ds.
\end{align*}
The Young inequality on the convolution yields
$$
\Vert AW(F_n) - AW(F)\Vert_{L^q(0,T;X)}
\le C\Vert t^{\alpha\ep-1}\Vert_{L^1(0,T)}
\left( \int^T_0 \Vert \AEAE (F_n-F)(s) \Vert^q ds\right)
^{\frac{1}{q}}
$$
$$
\le C\Vert F-F_n\Vert_{L^q(0,T;\DDD(\AEAE))}. 
$$
Therefore,
$$
\Vert AW(F_n) - AW(F)\Vert_{L^q(0,T;X)} \rrrr 0
\quad \mbox{as $n\to \infty$}.                       \eqno{(4.10)}
$$
Since $\pppa W(F_n) = AW(F_n)$ in $(0,T)$, we see
$$
\pppa W(F_n) \rrrr AW(F)\quad \mbox{in $L^q(0,T;X)$ as 
$n\to \infty$.}                                     \eqno{(4.11)}
$$
Similarly to (4.10), we can verify 
$$
\Vert W(F_n) - W(F)\Vert_{L^q(0,T;X)} \rrrr 0 \quad 
\mbox{as $n\to \infty$}.
$$
Since $\pppa$ is a closed operator with the domain $W_{\alpha,q}(0,T;X)
\subset L^q(0,T;X)$ to $L^q(0,T;X)$, it follows from (4.11) that  
$W(F) \in \DDD(\pppa)$ and 
$$
\pppa W(F) = \lim_{n\to\infty} \pppa W(F_n) = AW(F) \quad
\mbox{in $L^q(0,T;X)$}.
$$
Therefore, part (1) is proved for $a=0$ and $F \in L^q(0,T;\DDD(\AEAE))$.
Thus the proof of Theorem 4.1 is complete.
$\blacksquare$
\\
Now we proceed to 
\\
{\bf Proof of Proposition 4.1.}
\\
We will consider only the case $0<t<T$ and the arguments for
the cases $t=0$ and $t=T$ are the same. 
First we will prove the proposition for $a \in \DDD(A)$.

By (4.2) we have 
$$
G(t')a - G(t)a
= \frac{1}{2\pi i}\int_{\Gamma} (e^{\la t'} - e^{\la t})
\la^{\alpha-1}(\la^{\alpha} - A)^{-1}a d\la.
$$
We directly see
$$
\lim_{t'\to t}(e^{\la t'} - e^{\la t}) \la^{\alpha-1}(\la^{\alpha} - A)^{-1}a 
= 0 \quad \mbox{for all fixed $\la \in \Gamma$}.
$$
Moreover, by (2.7) we obtain
$$
\vert e^{\la t'}\vert, \,\, \vert e^{\la t}\vert
\le
\left\{ \begin{array}{rl}
& e^{-\vert \la\vert t\delta} + e^{-\vert \la\vert t'\delta}, 
\quad \la \in \Gamma_1 \cup \Gamma_3, \\
& e, \quad \la \in \Gamma_2.
\end{array}\right.
$$
Hence Condition ($\mathcal{A}$) yields 
$$
\Vert (e^{\la t'} - e^{\la t}) \la^{\alpha-1}(\la^{\alpha} - A)^{-1}a \Vert
\le
\left\{ \begin{array}{rl}
& Ce^{-\vert \la\vert\min\{t, t'\}\delta}\frac{1}{\vert \la\vert}, \quad
\la \in \Gamma_1 \cup \Gamma_3, \\
& \frac{C}{\vert \la\vert}, \quad \la \in \Gamma_2.
\end{array}\right.
                                            \eqno{(4.12)}
$$
By recalling the definition (2.3) of $\Gamma$, the right-hand side of (4.12) 
is in $L^1(\Gamma)$.  Therefore, the Lebesgue convergence theorem completes 
the proof of Proposition 4.1.
$\blacksquare$
\\
\vspace{0.2cm}
\\
{\bf 4.2. Uniqueness of solutions}
\\
The uniqueness of solution 
should be established within the class of solutions.
We discuss the uniqueness in several cases.

If we consider the solution for all $t>0$ with some growth rate 
as $t\to\infty$, then we can easily prove the uniquenss:
\\
{\bf Proposition 4.2.}
\\
{\it 
If $u\in W_{\alpha,q}(0,T;X) \cap L^q(0,T;\DDD(A))$ for all $T>0$ and
$\Vert u(t)\Vert_{X} = O(t^m)$ with some $m>0$ as $t\to\infty$
and $\pppa u(t) = Au(t)$ for $t>0$, then $u=0$ for $t>0$.
}
\\
{\bf Proof of Proposition 4.2.}
\\
By Lemma 2.4, we have
$$
L(\pppa u)(\la) = \la^{\alpha}Lu(\la) \quad \mbox{for all
$\la>0$.}
$$
Here we recall the Laplace transform $Lv(t) := \int^{\infty}_0
e^{-\la t}v(t) dt$.
Therefore, $\la^{\alpha}(Lu)(\la) = A(Lu)(\la)$ for all $\la>0$,
that is,
$$
(A-\la^{\alpha})(Lu)(\la) = 0 \quad \mbox{for all
$\la>0$.}
$$
Hence, $(Lu)(\la) = 0$ if $\la^{\alpha} \in \rho(A)$.

Since $\rho(A) \supset \Sigma_{\gamma}:= \{ z\in \C;\, 
\vert \mbox{arg}\, z\vert < \gamma,\, z\ne 0\}$ with some
$\gamma\in \left(\frac{\pi}{2}\, \pi\right)$ by 
{\it Condition ($\AAA$)}, if $\la \in \Sigma_{\gamma}$, then 
$\la^{\alpha} \in \Sigma_{\gamma} \subset \rho(A)$ by 
$0<\alpha<1$.  Therefore, $Lu(\la) = 0$ for all $\la \in 
\Sigma_{\gamma}$.
By $\{ \la > 0\} \subset \rho(A)$, we see that 
$(Lu)(\la) = 0$ for $\la > 0$.  The injectivity of the
Laplace transform yields $u(t) = 0$ for $t>0$.
Thus the proof of Proposition 4.2 is completed.
$\blacksquare$
\\

Proposition 4.2 implies that if
$$
\left\{ \begin{array}{rl}
& \pppa (u-a)(t) = Au(t) + F(t), \\
& \pppa (\www{u}-\www{a})(t) = A\www{u}(t) + \www{F}(t), \quad t>0,
\end{array}\right.
                                   \eqno{(4.13)}
$$
and $a=\www{a}$ and $F = \www{F}$ for $t>0$, then 
$u = \www{u}$ in $\OOO\times (0,\infty)$.
More precisely, Proposition 4.2 does not imply $u=\www{u}$ in 
$\OOO\times (0,T)$ from (4.13) only for $0<t<T$.

The case of $q=2$ and a Hilbert space $X$ makes the uniqueness easy.
\\
{\bf Proposition 4.3.}
\\
{\it
Let $X$ be a Hilbert space with scalar product $(\cdot,\cdot)_X$ and 
$$
q\ge 2                  
$$
and there exists a constant $C\ge 0$ such that 
$$
(Av, v)_X \le C\Vert v\Vert_X^2 \quad \mbox{for all $v \in \DDD(A)$}.
                                                      \eqno{(4.14)}
$$
If $u\in W_{\alpha,q}(0,T;X) \cap L^q(0,T;\DDD(A))$ satisfies 
$\pppa u = Au$ in $(0,T)$, then $u=0$ in $(0,T)$.
}

An elliptic opetator $A$ defined by (1.7) attached with 
the Dirichlet or the Neumann or the Robin boundary condition, 
satisfies (4.14).
\\
{\bf Proof of Proposition 4.3.}
\\
Since $W_{\alpha,q}(0,T;X) \cap L^q(0,T;\DDD(A))
\subset W_{\alpha,2}(0,T;X) \cap L^2(0,T;\DDD(A))$, it 
suffices to prove the uniqueness for $q=2$.

First we note the coercivity:
\\
{\bf Lemma 4.1.}
\\
{\it
Let $u\in W_{\alpha,2}(0,T;X)$.  Then
$$
\frac{1}{\Gamma(\alpha)}\int^t_0 (t-s)^{\alpha-1}
(\ppp_s^{\alpha} u(s),\, u(s))_X ds \ge \frac{1}{2}\Vert u(t)\Vert^2_X.
$$
}
As for the proof, see for example Theorem 3.4 (ii) in 
\cite{KRY}.
\\
{\bf Proof of Proposition 4.3.}
\\
By $\pppa u = Au$ and (4.14), in terms of Lemma 4.1,
we obtain
\begin{align*}
& 0 = \frac{1}{\Gamma(\alpha)}\int^t_0 (t-s)^{\alpha-1}
(\ppp_s^{\alpha} u(s), u(s))_Xds 
- \frac{1}{\Gamma(\alpha)}\int^t_0 (t-s)^{\alpha-1}
(Au(s), u(s))_X ds\\
\ge& \frac{1}{2}\Vert u(t)\Vert_X^2 
- C_1\int^t_0 (t-s)^{\alpha-1} \Vert u(s)\Vert^2_X ds,
\quad 0<t<T.
\end{align*}
The generalized Gronwall inequality (e.g., Henry \cite{H},
\cite{KRY}) yields $u=0$ in $(0,T)$, and so the proof of 
Proposition 4.3 is complete.
$\blacksquare$
\\

For the next result on the uniqueness, we consider $X = L^p(\OOO)$ with 
$1<p<\infty$ and an elliptic operator $A$ in $L^p(\OOO)$ defined by (1.7) and 
(1.8).  For simplicity, we further assume that all the coefficients 
$a_{k\ell}, b_k, c \in C^{\infty}(\ooo{\OOO})$ and $c\le 0$ in $\OOO$.
 
Then, we can prove
\\
{\bf Proposition 4.4.}
\\
{\it 
Let $X = L^p(\OOO)$ with $1<p<\infty$ and let 
$1 \le q \le \infty$, and let $A$ be defined by (4.16).
If $u \in W_{\alpha,q}(0,T;X) \cap L^q(0,T;\DDD(A))$ satisfies 
$\pppa u = Au$ in $(0,T)$, then $u=0$ in $\OOO\times (0,T)$.
}
\\
{\bf Proof of Proposition 4.4.}
\\
We know an estimate 
$$
\Vert v\Vert_{W^{2,p}(\OOO)} \le C(\Vert Av\Vert_{L^p(\OOO)}
+ \Vert v\Vert_{L^p(\OOO)}) \quad \mbox{for $v \in \DDD(A)$}
                      \eqno{(4.15)}
$$
(e,g., Ladyzhenskaya and Ural'tseva \cite{LU},
or Theorem 3.1 (p.212) in Pazy \cite{Pa}).
Taking into consideration that all the coefficients of 
$A$ are in $C^{\infty}(\ooo{\OOO})$, we can verify that
$A^{-\ell}a \in W^{2\ell,p}(\OOO)$ for each $\ell \in \N$ and 
$a\in L^p(\OOO)$.  Moreover, we can readily see that 
$J^mg \in W^{m,q}(0,T;L^p(\OOO))$ for each $m\in \N$ and 
$g \in L^q(0,T;L^p(\OOO))$.

In view of the Sobolev embedding,
for $1 < p < \infty$ and $1 \le q \le \infty$, we choose 
sufficiently large $\ell, m \in \N$ such that 
$$
W^{m,q}(0,T;W^{2\ell,p}(\OOO)) \subset W^{1,2}(0,T;L^2(\OOO)) \cap
L^2(0,T;H^2(\OOO)\cap H^1_0(\OOO)).
$$
Then,
$$
w:= J^mA^{-\ell}u \in 
W^{1,2}(0,T;L^2(\OOO)) \cap L^2(0,T;H^2(\OOO)\cap H^1_0(\OOO))
                                                  \eqno{(4.16)}
$$
for $u \in L^q(0,T;L^p(\OOO))$.

Noting that $J^m\pppa u = \pppa J^mu$ for 
$u \in \DDD(\pppa) = W_{\alpha,q}(0,T;L^p(\OOO))$, we operate
$J^mA^{-\ell}$ to the equation $\pppa u = Au$ in $L^q(0,T;L^p(\OOO))$,
so that we obtain
$$
\pppa w = Aw \quad \mbox{in $L^2(0,T;L^2(\OOO))$.}     \eqno{(4.17)}
$$
In terms of (4.16), for (4.17) we can apply the uniqueness of solution 
$w \in L^2(0,T;\DDD(A)) \cap W_{\alpha,2}(0,T;L^2(\OOO))$
(e.g., \cite{KRY}, \cite{Za}), and so $w=J^mA^{-\ell} u =0$ in 
$\OOO\times (0,T)$.
As is readily verified, the operator $J^m$ and $A^{-\ell}$ are injective, and
we reach $u=0$ in $\OOO\times (0,T)$.
Thus the proof of Proposition 4.4 is complete.
$\blacksquare$
\\

In each solutions space in Propositions 4.2 - 4.4, there 
exists at most one strong solution to (4.3) and 
it is given by (4.7).
\\

Theorem 4.1 in the case of $q=2$ and $X=L^2(\OOO)$ is 
the same as many previous results (e.g., Jin \cite{J},
Kubica, Ryszewska and Yamamoto \cite{KRY}, Sakamoto and 
Yamamoto \cite{SY}, for example), and the existing works
are mainly based on the eigenfunction expansions of the 
solutions.

Our approach produces solution formula (4.7), which is 
convenient for applications also in non-Hilbert space $X$, 
as are shown in the next sections.
We conclude this section with the unique existence of 
$u$ to (4.1) omitting the assumption $0 \in \rho(A)$.
\\
{\bf Theorem 4.3.}
\\
{\it
We assume that there exists a constant $C_0>0$ such that 
$$
\mbox{$A_0:= A-C_0$ with $\DDD(A_0) = \DDD(A)$ satisfies 
{\it Condition ($\AAA$)} but
not necessarily $0 \in \rho(A)$.}                     \eqno{(4.18)}
$$
Moreover let (4.4) and (4.5) hold.
Then, for $a \in \DDD((-A_0)^{\mu})$ and $F \in L^q(0,T; \DDD((-A_0)^{\ep}))$,
the problem (4.3) possesses a strong solution 
$u \in L^q(0,T;\DDD(A_0))$ satisfying $u-a \in W_{\alpha,q}(0,T;X)$ and
$$
\Vert u\Vert_{L^q(0,T;\DDD(A_0))} + \Vert u-a\Vert_{W_{\alpha,q}(0,T;X)}
\le C(\Vert a\Vert_{\DDD((-A_0)^{\mu})} + \Vert F \Vert
_{L^q(0,T; \DDD((-A_0)^{\ep}))}).                \eqno{(4.19)}
$$
}
\\
{\bf Remark.}
Unlike the case $\alpha=1$, we cannot shift the zeroth-order term
by a simple transform $u \mapsto ue^{-C_0t}$.
Thus for the case $0<\alpha<1$, we have to resort to the iteration 
method.
\\
{\bf Proof of Theorem 4.3.}
\\
We rewrite (4.1) as
$$
\pppa (u(t)-a) = A_0u(t) + C_0u(t) + F(t), \quad 0<t<T.
$$
We can construct the operator $G(t)$ and $K(t)$ by (4.2) and (4.3) for $A_0$.  
We introduce an iteration scheme:
$$
\left\{ \begin{array}{rl}
& u_0(t) = G(t)a + \int^t_0 K(t-s)F(s) ds, \cr \\
& u_{n+1}(t) = G(t)a + \int^t_0 C_0K(t-s)u_n(s) ds
+ \int^t_0 K(t-s)F(s) ds, \quad n\in \N \cup \{0\}.
\end{array}\right.
                                  \eqno{(4.20)}
$$
By Theorem 4.1 (2), we see 
$$
u_0 \in L^q(0,T;\DDD(A_0)), \quad u-a \in W_{\alpha,q}(0,T;X)
                                            \eqno{(4.21)}
$$
and
$$
\Vert A_0u_0(t)\Vert \le Ct^{-\alpha(1-\mu)}\Vert (-A_0)^{-\mu}a\Vert
+ C\int^t_0 (t-s)^{\alpha\ep-1} \Vert (-A_0)^{\ep} F(s)\Vert ds,
\quad 0<t<T.
$$
By (4.18), we have
\begin{align*}
& \Vert A_0(u_2(t) - u_1(t))\Vert
= \left\Vert \int^t_0 A_0C_0K(t-s) (u_1(s) - u_0(s))ds \right\Vert\\
= &\left\Vert \int^t_0 C_0K(t-s) A_0(u_1(s) - u_0(s))ds \right\Vert\\
\le& C\int^t_0 \Vert K(t-s)\Vert \Vert A_0(u_1(s) - u_0(s))\Vert ds 
\le C\int^t_0 (t-s)^{\alpha-1}\Vert A_0(u_1(s) - u_0(s))\Vert ds.
\end{align*}
Hence, repeatedly applying similar estimates, we have
\begin{align*}
& \Vert A_0(u_3(t) - u_2(t))\Vert
\le C\int^t_0 (t-s)^{\alpha-1}
C\left(\int^s_0 (s-\xi)^{\alpha-1} \Vert A_0(u_1(\xi) - u_0(\xi))\Vert 
d\xi \right)ds\\
=& C^2\int^t_0 \left( \int^t_{\xi} (t-s)^{\alpha-1}
 (s-\xi)^{\alpha-1} ds \right) \Vert A_0(u_1(\xi) - u_0(\xi))\Vert
d\xi\\
=& C^2\int^t_0 \frac{\Gamma(\alpha)^2}{\Gamma(2\alpha)}
 (t-s)^{2\alpha-1} \Vert A_0(u_1(\xi) - u_0(\xi))\Vert d\xi.
\end{align*}
Here we exchanged the orders of the integrals and 
$$
\int^t_{\xi} (t-s)^{\alpha-1}(s-\xi)^{\alpha-1} ds
= \int^{t-\xi}_0 (t-\xi-\eta)^{\alpha-1}\eta^{\alpha-1} d\eta,
$$
which is verified by the change of the variables $\eta = s-\xi$.

Hence, 
$$
\Vert A_0(u_3(t) - u_2(t))\Vert
\le \frac{(C\Gamma(\alpha)^2}{\Gamma(2\alpha)}
\int^t_0  (t-s)^{2\alpha-1} \Vert A_0(u_1(s) - u_0(s))\Vert ds.
$$
Continuing the estimation, we can reach 
$$
\Vert A_0(u_{n+1}(t) - u_n(t))\Vert
\le \frac{(C\Gamma(\alpha)^n}{\Gamma(n\alpha)}
\int^t_0  (t-s)^{n\alpha-1} \Vert A_0(u_1(s) - u_0(s))\Vert ds, \quad
0<t<T.
$$
Therefore, in terms of (4.21), the Young inequality on the 
convolution yields
\begin{align*}
& \Vert A_0(u_{n+1} - u_n)\Vert_{L^q(0,T;X)}
= \Vert u_{n+1} - u_n\Vert_{L^q(0,T;\DDD(A_0))} \\
\le & \frac{(C\Gamma(\alpha))^n}{\Gamma(n\alpha)}
\frac{T^{n\alpha}}{n\alpha} \Vert A_0(u_1 - u_0)\Vert_{L^q(0,T;X)}
\end{align*}
$$
\le \frac{(C\Gamma(\alpha)T^{\alpha})^n}{\Gamma(n\alpha+1)}
(\Vert A_0u_1\Vert_{L^q(0,T;\DDD(A_0))} + \Vert u_0\Vert
_{L^q(0,T;\DDD(A_0))})
                                       \eqno{(4.22)}
$$
for each $n\in \N$.
Since 
$$
\lim_{n\to \infty} 
\frac{(C\Gamma(\alpha)T^{\alpha})^{n+1}}{\Gamma((n+1)\alpha+1)}
\left( \frac{(C\Gamma(\alpha)T^{\alpha})^n}{\Gamma(n\alpha+1)}
\right)^{-1}
= C\Gamma(\alpha)T^{\alpha} \lim_{n\to \infty}
\frac{\Gamma(n\alpha+1)}{\Gamma((n+1)\alpha+1)} = 0,
$$
we see 
$$
\sum_{n=0}^{\infty} 
\frac{(C\Gamma(\alpha)T^{\alpha})^n}{\Gamma(n\alpha+1)} < \infty.
$$
It follows from (4.22) that 
$$
\left\{ \begin{array}{rl}
& \Vert A_0u_n\Vert_{L^q(0,T;X)} 
\le C(\Vert A_0u_1\Vert_{L^q(0,T;X)} + \Vert A_0u_0\Vert
_{L^q(0,T;X)}) \quad \mbox{for each $n\in \N$}, \cr\\
& A_0u_n \rrrr A_0u \quad \mbox{in $L^q(0,T;X)$ as $n\to \infty$}.
\end{array}\right.
                                    \eqno{(4.23)}
$$
By Theorem 4.1, we have
$$
\pppa (u_{n+1} - a) = A_0u_{n+1} + C_0u_n + F(t),\quad 0<t<T.
$$
In view of (4.3), we see that 
$\lim_{n\to \infty} \pppa (u_{n+1} - a)=: \www{u}$ is convergent 
in $L^q(0,T;X)$ and $\www{u} = A_0u + C_0u + F(t)$ for $0<t<T$.
By the closedness of the operator $\pppa: W_{\alpha,q}(0,T;X) 
\subset L^q(0,T;X) \rrrr L^q(0,T;X)$, we obtain
$$
\pppa (u-a) = A_0u + C_0u + F(t), \quad 0<t<T.
$$
With (4.23), the proof of Theorem 4.3 is complete.
$\blacksquare$
\\
\section{Semilinear equations}

Thanks to solution formula (4.7) with estimates 
(4.8) and (4.9), we can apply a perturbation method by which 
we can discuss semilinear equations.  Such a treatment is quite 
standard for the case $\alpha=1$, and we can refer to 
tremendously many works.  Here we refer only to 
Henry \cite{H}, Pazy \cite{Pa}, Yagi \cite{Ya}.

In order to highlight the essence, we consider a simple but typical 
case: let $1< p < \infty$, $X=L^p(\OOO)$ and 
$$
\AAA v(x) = \Delta v(x) = \sum_{k=1}^d \ppp_k^2v(x), \quad
\DDD(A) = \{ v\in W^{2,p}(\OOO);\, v\vert_{\ppp\OOO} = 0\}.
$$
Here $\OOO \subset \R^d$ is a bounded smooth domain.

We consider 
$$
\pppa (u(t) - a) = Au(t) + F(u(t)), \quad 0<t<T.    \eqno{(5.1)}
$$
Here $F$ is assumed to satisfy the conditions as stated later.
Our main concern is the unique existence of the solution to (5.1) 
with small $T>0$.

It is known (e.g., Pazy \cite{Pa}, Tanabe \cite{Ta}) that for
$0 < \gamma \le 1$, the fractional power $(-A)^{\gamma}$ of $-A$ 
is well-defined.
Henceforth, we choose $\gamma, \gamma' \in (0,1)$ such that 
$$
\gamma' p > \frac{d}{2}, \quad 0<\gamma' < \gamma.              \eqno{(5.2)}
$$
Then, we have 
$$
\DDD((-A)^{\gamma}) \subset W^{2\gamma',p}(\OOO)
$$
and there exists a constant $C = C(\gamma, \gamma') > 0$ such that 
$$
\Vert v\Vert_{W^{2\gamma',p}(\OOO)} \le C\Vert (-A)^{\gamma}v\Vert
_{L^p(\OOO)} \quad \mbox{for $v \in \DDD((-A)^{\gamma})$}
                                                           \eqno{(5.3)}
$$
(e.g., Theorem 1.6.1 (p.39) in \cite{H}).
By (5.2) and the Sobolev embedding, we see
$$
\DDD((-A)^{\gamma}) \subset W^{2\gamma',p}(\OOO) \subset C(\ooo{\OOO}).
                                                  \eqno{(5.4)}
$$

Now we introduce the conditions posed on the semilinear term 
$F$ of equation (5.1). 
For a constant $M>0$, there exists a constant $C = C(M) >0$ such 
that
$$
\Vert F(v)\Vert \le C(M), \quad 
\Vert F(v_1) - F(v_2)\Vert \le C(M)\Vert v_1-v_2\Vert_{\DDD((-A)^{\gamma})}
$$
$$
\mbox{if $\Vert v\Vert_{\DDD((-A)^{\gamma})}$, 
$\Vert v_1\Vert_{\DDD((-A)^{\gamma})}$,  
$\Vert v_2\Vert_{\DDD((-A)^{\gamma})}\le M$}.
                                         \eqno{(5.5)}
$$
Henceforth, by $C>0$, $C_0, C_1 > 0$, etc., we denote generic constants, which 
are independent of the functions $u, v$, etc. under consideration, and 
we write 
$C_f$, $C(m)$ in the case where 
we need to specify a dependence on related quantities.
Before we state the main result of this section, we show an 
example of semilinear functions, which satisfies condition (5.5).
\\
{\bf Example.}\\
For $F\in C^1(\R)$, by setting $F(v):= F(v(x))$ for $x \in \OOO$,
we can define $F: \DDD(\DDDg)\, \longrightarrow L^p(\OOO)$ and 
$F$ satisfies (5.5) if
$$
\gamma p > \frac{d}{2} \quad \mbox{and} \quad 0<\gamma<1.   \eqno{(5.6)}
$$

Indeed, by (5.6) we can choose $0 < \gamma' < \gamma$ such that (5.2)
is satisfied.  Then,
$\Vert v\Vert_{\DDD(\DDDg)}\le M$ yields
$\vert v(x)\vert \le C_0M$ for $x \in \OOO$, so that 
$$
\Vert F(v)\Vert_{L^{\infty}(\OOO)} \le \max_{\vert \eta\vert\le C_0M}
\vert F(\eta)\vert =: C(M)
$$
for $\Vert v\Vert_{\DDD((-A)^{\gamma})} \le M$,
Therefore, the first condition in (5.5) holds true.
By using the mean value theorem, the second condition in (5.5) is 
satisfied.
$\blacksquare$
\\

Now we define a class of solutions:
\\
{\bf Definition 5.1 (mild solution)}.
\\
{\it 
We assume that $F(v)$ is well-defined for $v \in \DDD((-A)^{\gamma})$ 
with a constant $\gamma > 0$ satisfying (5.2).
We call $u=u(t) \in C([0,T]; \DDD(\DDDg))$ a mild solution to (5.1) if
$$
u(t) = G(t)a + \int^t_0 K(t-s)F(u(s)) ds, \quad 0<t<T.   \eqno{(5.7)}
$$
}

The notion of mild solution is standard in the case of $\alpha=1$,
we can refer for example to Definition 1.1 (p.184) in Pazy \cite{Pa}.
We remark that sufficiently smooth $u$ satisfies (5.1) for 
$(x,t) \in \OOO\times (0,T)$, then it is a mild solution, and moreover 
we can prove that a sufficiently smooth mild solution $u$ satisfies
(5.1) for all $(x,t)\in \OOO\times (0,T)$.
For (5.1), we can study which smoothness condition on $a$ and $F$
allows a mild solution to 
possesses the regularity such as $\pppa (u-a), \, Au \in C([0,T];X)$, which 
means that $u$ is a classical solution.  However, we here do not discuss the 
details.

In a special case $p=2$, condition (5.2) requires 
$d < 4\gamma' < 4$, and for $d\ge 4$, 
we cannot treat higher spatial dimensions directly as long as 
we remain in the $L^2(\OOO)$-space.
Within our framework, by choosing $p>2$, we can uniformly treat
the higher dimensional case $d > 4$, and this is a main advantage of the 
$L^p$-approach.
\\

Now we state our main result regarding the local unique 
existence of a mild solution to an initial boundary value problem for 
semilinear fractional diffusion equation (5.1). 
\\
{\bf Theorem 5.1.}
\\
{\it
Let a semilinear term $F$ satisfy conditions (5.5) with $M>0$ and 
$\Vert a\Vert_{\DDD(\DDDg)} \le M$.  
Then there exists a constant $T=T(M)>0$
such that (5.1) possesses a unique mild solution 
$$
u_a \in C([0,T]; \DDD(\DDDg)).                \eqno{(5.8)}
$$
Moreover, there exists a constant $C(M)>0$, such that
$$
\Vert u_a - u_b\Vert_{C([0,T];\DDD(\DDDg))}
\le C\Vert a-b\Vert_{\DDD(\DDDg)}               \eqno{(5.9)}                  
$$
provided that $\Vert a\Vert_{\DDD(\DDDg)}, \Vert b\Vert_{\DDD(\DDDg)} 
\le M$.
}

Theorem 5.1 is similar to the ones 
well-known for the partial differential equations of parabolic type 
that correspond to the case $\alpha=1$ 
(see, e.g.,  \cite{H} or \cite{Pa}). 
\\
\vspace{0.2cm}
\\
{\bf Proof of Theorem 5.1.}
\\
{\bf First Step.}
We prove the uniqueness of mild solution.
Let $a=0$ and let $u_k$, $k=1,2$ be mild solutions to 
$$
u_k(t) = \int^t_0 K(t-s) F(u_k(s)) ds, \quad k=1,2, \quad
0<t<T.
$$
Then, 
$$
\DDDg u_k(t) = \int^t_0 \DDDg K(t-s) F(u_k(s)) ds, \quad k=1,2, \quad
0<t<T.
$$
By (4.9) and (5.5), we have
\begin{align*}
& \Vert \DDDg (u_1(t)-u_2(t)) \Vert 
\le C\int^t_0 (t-s)^{\alpha(1-\gamma)-1} 
\Vert F(u_1(s)) - F(u_2(s))\Vert ds\\
\le& C_M\int^t_0 (t-s)^{\alpha(1-\gamma)-1} 
\Vert \DDDg(u_1(s) - u_2(s))\Vert ds, \quad 0<t<T,
\end{align*}
where we set $M:= \max\{ \Vert u_1\Vert_{C([0,T];\DDD(\DDDg))}, \,
\Vert u_2\Vert_{C([0,T];\DDD(\DDDg))} \}$.
The generalized Gronwall inequality implies 
$\Vert \DDDg( u_1(t) - u_2(t))\Vert = 0$ for $0<t<T$, that is, 
$u_1 = u_2$ in $\OOO\times (0,T)$.  Thus the uniqueness of mild solution
is proved.
$\blacksquare$
\\
{\bf Second Step.}
In this step, we prove
\\
{\bf Lemma 5.1.}
\\
{\it Let $H \in C([0,T]; L^p(\OOO))$.  Then
$$
\int^t_0 \DDDg K(t-s)H(s) ds \in C([0,T];L^p(\OOO)),
$$
that is,
$$
\int^t_0 K(t-s)H(s) ds \in C([0,T];\mathcal{D}(\DDDg)).
$$
}
\\
{\bf Proof of Lemma 5.1.}
\\
Let $0< \eta < t \le T$.  We have the representation
\begin{align*}
& \int^t_0 \DDDg K(t-s)H(s) ds - \int^{\eta}_0 \DDDg K(\eta-s)
H(s) ds\\
=& \int^t_0 \DDDg K(s)H(t-s) ds - \int^{\eta}_0 \DDDg K(s)
H(\eta-s) ds\\
=& \int^t_{\eta} \DDDg K(s)H(t-s) ds 
+ \int^{\eta}_0 \DDDg K(s) (H(t-s) - H(\eta-s)) ds\\
=: &I_1 + I_2.
\end{align*}
For the first integral, by (4.9) in Theorem 4.1 and $\gamma < 1$, 
we have the relations
\begin{align*}
&\Vert I_1\Vert
\le C\int^t_{\eta} s^{\alpha(1-\gamma)-1} \max_{0\le s\le t}
\Vert H(t-s)\Vert ds\\
\le& C\Vert H\Vert_{C([0,T];L^p(\OOO))}
\frac{t^{\alpha(1-\gamma)} - \eta^{\alpha(1-\gamma)}}
{\alpha(1-\gamma)} \, \longrightarrow \, 0
\end{align*}
as $\eta \uparrow t$.  

Next, in terms of (4.9), we similarly obtain
\begin{align*}
& \Vert I_2\Vert
=\left\Vert \int^{\eta}_0 \DDDg K(s) (H(t-s) - H(\eta-s)) ds\right\Vert\\
\le & C\int^{\eta}_0 s^{(1-\gamma)\alpha-1} 
\Vert H(t-s) - H(\eta-s)\Vert ds.
\end{align*}
For $H\in C([0,T];L^p(\OOO))$, the function
$$
\vert s^{(1-\gamma)\alpha-1} \vert
\Vert H(t-s) - H(\eta-s)\Vert 
\le s^{(1-\gamma)\alpha-1} \times 2 \Vert H\Vert
_{C([0,T];L^p(\OOO))}
$$
is an integrable function with respect to $s \in (0,\eta)$ and 
$$
\lim_{\eta \uparrow t} s^{(1-\gamma)\alpha-1} 
\Vert H(t-s) - H(\eta-s)\Vert = 0
$$
for almost all $s\in (0,\eta)$.
Hence, the Lebesgue convergence theorem implies  
$\lim_{\eta \uparrow t} \Vert I_2\Vert = 0$ and the proof of 
Lemma 5.1 is completed.
$\blacksquare$ 
\\
{\bf Third Step.}
We proceed to completion of the proof of Theorem 5.1.
In view of Proposition 4.1, the inclusion $a \in \DDD(\DDDg)$ 
implies
$$
G(t)a \in C([0,T];\DDD(\DDDg)).                     \eqno{(5.10)}
$$                       
Indeed, 
$$
 \Vert \DDDg (G(t')a - G(t)a)\Vert
= \Vert G(t')(\DDDg a) - G(t)\DDDg a)\Vert,
$$
so that $\lim_{t'\to t} G(t')(\DDDg a) = G(t)(\DDDg a)$,
which implies (5.10).
\\

We set
$$
Lv(t) = G(t)a + \int^t_0 K(t-s)F(v(s)) ds \quad \mbox{for 
$v \in C([0,T];\DDD(\DDDg))$},
$$
and
$$
V:= \{ v\in C([0,T];\DDD(\DDDg));\, 
\Vert v - Ga\Vert_{C([0,T];\DDD(\DDDg))} \le M
                                               \eqno{(5.11)}
$$
with some constant $M>0$.

By (5.4) and (5.5), we obtain $F(v) \in C([0,T];L^p(\OOO))$ for
$v \in C([0,T];\DDD(\DDDg))$. Therefore, applying Lemma 5.1 where we choose
$H(s):= F(v(s))$, in view of (5.10), we reach 
$$
Lv \in C([0,T]; \DDD((-A)^{\gamma})) \quad 
\mbox{for $v \in C([0,T]; \DDD((-A)^{\gamma}))$}.
$$
Therefore, $L C([0,T]; \DDD((-A)^{\gamma})) \subset 
C([0,T]; \DDD((-A)^{\gamma}))$.

Choosing $T>0$ sufficiently small, we will prove
\\
(i) $LV \subset V$, where $V$ is the set defined by (5.11).
\\
(ii) There exists a constant $\rho \in (0,1)$ such that  
$$
\Vert Lu_1 - Lu_2\Vert_{C[0,T];\DDD(\DDDg))} 
\le \rho\Vert u_1 - u_2\Vert_{C[0,T];\DDD(\DDDg))}, \quad 0<t<T
\quad \mbox{for any $u_1, u_2 \in V$}.
$$
\\
{\bf Proof of (i).}
Let $u \in V$.  Now we consider the expression 
$$
\DDDg (Lu(t) - G(t)a)
= \int^t_0 \DDDg K(t-s) F(u(s)) ds.     \eqno{(5.12)}
$$
For any $u \in V$, using the norm estimates 
$$
\Vert a\Vert_{\mathcal{D}(\DDDg)}
= \Vert \DDDg a\Vert \le M, \quad
\Vert u-G(\cdot)a\Vert_{C([0,T];\mathcal{D}(\DDDg))} \le M,
$$
by Proposition 2.1 (ii) we obtain
$$
\Vert u(t)\Vert_{\DDD(\DDDg)} \le M + \Vert \DDDg G(t)a\Vert 
= M + \Vert G(t)\DDDg a\Vert 
\le M + C_1M =: C_2M                       \eqno{(5.13)}
$$
for each $v \in V$.
The first condition from (5.5) implies that
$$
\Vert F(u(t)) \Vert \le C(C_2M) \quad \mbox{for all $u\in V$
and $0<t<T$.}                                                  \eqno{(5.14)}
$$
Applying (5.14) to (5.12), by means of (4.9),
we obtain the norm estimates
\begin{align*}
& \Vert Lu(t) - G(t)a \Vert_{\DDD(\DDDg)} 
= \left\Vert \int^t_0 \DDDg K(t-s)F(u(s)) ds\right\Vert\\
\le & C\int^t_0 (t-s)^{(1-\gamma)\alpha-1} C(C_2M)ds
\le C_3\frac{t^{(1-\gamma)\alpha}}{(1-\gamma)\alpha}
\le C_3\frac{T^{(1-\gamma)\alpha}}{(1-\gamma)\alpha}.
\end{align*}
The constant $C_3>0$ depends on $M>0$ but is independent 
on $T>0$.
Therefore, choosing $T>0$ sufficiently small, we complete the proof of 
the property (i).  $\square$
\\
{\bf Proof of (ii).}
Estimate (5.13) yields that $\Vert u_1(t)\Vert_{\DDD(\DDDg)} \le C_2M$ and 
$\Vert u_2(t)\Vert_{\DDD(\DDDg)} \le C_2M$ for any $u_1, u_2 \in V$.
Condition (5.5) leads to 
$$
\Vert F(u_1(s)) - F(u_2(s))\Vert \le C(C_2M)\Vert u_1(s)-u_2(s)\Vert
_{\DDD(\DDDg)}, \quad 0<s<T.
$$
Hence, we have
\begin{align*}
& \Vert Lu_1(t) - Lu_2(t)\Vert_{\DDD(\DDDg)} 
= \left\Vert  \int^t_0 \DDDg K(t-s)(F(u_1(s)) - F(u_2(s))) ds\right\Vert\\
\le& C(C_2M)\int^t_0 (t-s)^{\alpha(1-\gamma)-1}\Vert (u_1-u_2)(s)\Vert
_{\DDD(\DDDg)} ds\\
\le & C_4T^{\alpha(1-\gamma)}
\sup_{0<s<T} \Vert u_1(s) - u_2(s)\Vert_{\DDD(\DDDg)}.
\end{align*}
In the last inequality, the constant $C_4$ is independent of $T$, and thus 
we can  choose a sufficiently small constant $T>0$ satisfying 
the inequality
$$
\rho:= C_4T^{\alpha(1-\gamma)} < 1.
$$
The proof of the property (ii) is completed.
$\blacksquare$

Thus the contraction theorem for small $T>0$ completes the proof of 
Theorem 5.1.
$\blacksquare$
\section{Weak solution and smoother solution}

{\bf 6.1. Weak solution}

We return to the linear equation and construct a weak solution.
For it, the estimate (4.9) is helpful.  For concentrating on the 
essence, we discuss the case where the initial value is zero, that is,
$a=0$:
$$
\pppa u(t) = Au(t) + F(t), \quad 0<t<T \quad \mbox{in $X$}.    \eqno{(6.1)}
$$
Henceforth, for a Banach space $X$ over $\R$, let $X'$ be the 
dual space of $X$: the Banach space of all the bounded linear functionals
from $X$ to $\R$.  By $\, _{X'}\langle f, u \rangle_X 
= \langle f, u\rangle$ we denote the value of $f\in X'$ at $u\in X$.
Moreover let $A'$ be the dual operator to $A$, which is defined in $X'$
with values in $X'$.

Henceforth $(\pppa)'$ denotes the dual operator to 
$\pppa : W_{\alpha,q}(0,T;X) \longrightarrow L^q(0,T;X)$.

For $1 \le r \le \infty$, we define an operator $\tau: L^r(0,T;L^q(X) 
\longrightarrow L^r(0,T;X)$ by 
$(\tau v)(x,t) = v(x,T-t)$.
Then, we can prove
\\
{\bf Lemma 6.1.}
\\
{\it 
Let $1 < q < \infty$ and $\frac{1}{q} + \frac{1}{q'} = 1$.  Then, 
$$
\DDD((\pppa)') = \tau W_{\alpha,q'}(0,T;X'), \quad
(\pppa)'u = \tau\pppa (\tau u).
$$
}
For the proof, see e.g., Yamamoto \cite{Y18}, \cite{Y22}.

By the definition of the dual operator $(\pppa)'$, we see
$$
\langle u, \, (\pppa)'v \rangle 
= \langle \pppa u,\, v \rangle \quad \mbox{for $u\in \DDD(\pppa)$ 
and $v \in \DDD((\pppa)')$}.
$$

We define the space of test functions by
$$
\Psi:= \{ \psi\in C^{\infty}([0,T]; \DDD(A'))\,; \,\,  
\psi(\cdot,T) = 0 \quad \mbox{in $\OOO$} \}.
$$
Then we can readily prove 
$$
\Psi \subset \DDD((\pppa)').            \eqno{(6.2)}
$$
\\
{\bf Definition 6.1 (weak solution)}
\\
{\it
For $F\in L^q(0,T;X)$, we call $u$ a weak solution to (6.1), 
if $u\in L^q(0,T;X)$  and
$$
 \, _{L^q(0,T;X)}\langle u, \, ((\pppa)' - A')\psi \rangle_{L^{q'}(0,T;X')} 
= \, _{L^q(0,T;X)}\langle F,\, \psi \rangle_{L^{q'}(0,T;X')}
$$
for all $\psi \in \Psi$.
}

Henceforth, we write 
$$
\langle u,\, v\rangle
= \, _{L^q(0,T;X)}\langle u, \, v\rangle_{L^{q'}(0,T;X')}.
$$

We are ready to state the main result in Subsection 6.1.
\\
{\bf Theorem 6.1.}
\\
{\it 
Let $F \in L^q(0,T;X)$ for $1 \le q\le \infty$.  Then there exists a unique 
weak solution $u$ to (6.1) and 
$$
u(t) = \int^t_0 K(t-s)F(s) ds, \quad 0<t<T.
$$
}
\\
{\bf Proof of Uniqueness.}
\\
Let $u \in L^q(0,T;X)$ satisfy 
$$
\langle u, \, (\pppa)'\psi - A'\psi \rangle = 0 
\quad \mbox{for all $\psi \in \Psi$.}              \eqno{(6.3)}
$$
Then, we have to prove that $u=0$ in $\OOO\times (0,T)$.

We choose $\ell, m \in \N$ satisfying (4.16), that is,
$$
w := J^mA^{-\ell}u \in W_{1,2}(0,T;X) \cap L^2(0,T;H^2(\OOO)\cap H^1_0(\OOO)).
                                     \eqno{(6.4)}
$$

We first show two lemmata.
\\
{\bf Lemma 6.2.}
\\
{\it
$$ 
(A^{-\ell})'(J^m)' \psi \in \Psi \quad \mbox{for each 
$\psi \in \Psi$}.
$$
}
\\
{\bf Proof of Lemma 6.2.}
\\
It suffices to prove $(A^{-\ell})'v \in C^{\infty}([0,T];\DDD(A'))$ and
$(J^m)'v \in C^{\infty}([0,T];\DDD(A'))$ for 
$v \in C^{\infty}([0,T];\DDD(A'))$ and $v(T) = 0$. The verification of the
latter is straightforward. Indeed, as is directly shown,
$$
(J^m)'v(t) = \frac{1}{\Gamma(m)}\int^T_t (\xi-t)^{m-1}v(\xi) d\xi,
\quad 0<t<T
$$
which implies the latter inclusion.

As for the former inclusion, we remark that $(A^{-\ell})' = (A')^{-\ell}$.
If $y \in \DDD(A')$, then $y = (A')^{-1}z$ with some $z \in X'$, and so
$$
(A^{-\ell})'y = (A')^{-\ell}y = (A')^{-1}((A')^{-\ell}z) \in \DDD(A'),
$$
that is, $y \in \DDD(A')$ implies $(A^{-\ell})'y \in \DDD(A')$.
Thus the proof of Lemma 6.2 is complete.
$\blacksquare$
\\
{\bf Lemma 6.3.}
\\
$$ 
(J^m)'(\pppa)'\psi = (\pppa)'(J^m)'\psi \quad \mbox{for each 
$\psi \in \Psi$}.
$$
\\
{\bf Proof of Lemma 6.3.}
\\
We can directly see
$(J^m)'\psi(t) = (\tau J^m(\tau \psi)(t)$, where 
we recall $(\tau\psi)(t): = \psi(T-t)$ for $0<t<T$.

Therefore, since by the definition of $\pppa = (J^{\alpha})^{-1}$, we
can readily see that $J^m\pppa = \pppa J^m$ on $\DDD(\pppa)$, so that 
$$
(J^m)'((\pppa)'\psi(t)) = (J^m)'(\tau \pppa (\tau\psi))(t)
= \tau J^m(\tau \tau \pppa (\tau\psi))(t) 
= \tau(J^m\pppa (\tau \psi))(t).
$$
Consequently, since $\tau^2$ is the identity mapping, we derive 
$$
(J^m)'(\pppa)'\psi = \tau J^m ((\pppa(\tau \psi)))
= \tau \pppa J^m(\tau\psi)
= \tau\pppa \tau(\tau J^m(\tau\psi)) = (\pppa)'(J^m)'\psi.
$$
Thus the proof of Lemma 6.3 is complete.
$\blacksquare$ 

We return to the proof of the uniqueness.  We obtain 
$$
 \langle w,\, (\pppa)'\psi\rangle = \langle J^mA^{-\ell}u,\, (\pppa)'\psi
\rangle =  \langle u,\, (A^{-\ell})'(J^m)'(\pppa)'\psi \rangle
$$
and
$$
 \langle w,\, A'\psi\rangle = \langle u, \, (A^{-\ell})'(J^m)'A'\psi \rangle.
$$
Therefore, using Lemma 6.3, we have 
$$
 \langle w,\, ((\pppa)' - A')\psi\rangle 
= \langle u, \, ((\pppa)' - A')((A^{-\ell})'(J^m)'\psi)\rangle.
                                                  \eqno{(6.5)}
$$
In view of Lemma 6.2, we know 
$(A^{-\ell})'(J^m)'\Psi \subset \Psi$, which yields
$$
\langle u, \, ((\pppa)' - A')((A^{-\ell})'(J^m)'\psi)\rangle
  = 0
$$
for all $\psi\in \Psi$.
The application of (6.5) yields 
$$
\langle w, \, ((\pppa)' - A')\psi)\rangle = 0 \quad \mbox{for all 
$\psi \in \Psi$}.
$$
By (6.4), choosing $\psi \in C^{\infty}_0(\OOO \times (0,T)) \subset \Phi$, we derive 
$$
\int^T_0 \int_{\OOO} w(x,t) ((\pppa)' - A')\psi(x, t) dxdt = 0
$$
for all $\psi\in C^{\infty}_0(\OOO \times (0,T))$, which implies 
$$
(\pppa - A)w = 0                  \eqno{(6.6)}
$$
in the sense of distribution.  
In terms of (6.4), it follows that (6.6) holds in $L^2(0,T;L^2(\OOO))$.
Applying the uniqueness of solution in $W_{1,2}(0,T;X) 
\cap L^2(0,T;H^2(\OOO) \cup H^1_0(\OOO))$ (e.g., \cite{KRY}) to (6.6), we 
reach $w=0$ in $\OOO\times (0,T)$.  Since $J^mA^{-\ell}$ is injective, 
we obtain $u=0$ in $\OOO\times (0,T)$.  Thus the proof of the uniqueness
is complete.
$\blacksquare$
\\
{\bf Proof of Existence.}
\\
Henceforth we set 
$$
w(F)(t) = \int^t_0 K(t-s) F(s) ds, \quad 0<t<T.
$$
Then, 
\\
{\bf Lemma 6.4}
{\it
$$
\Vert w(F)\Vert_{L^q(0,T;X)} \le C\Vert F\Vert_{L^q(0,T;X)}
$$
for each $F \in L^q(0,T;X)$.
}
\\
{\bf Proof of Lemma 6.4.}
\\
Estimate (4.9) in Theorem 4.1 yields 
$$
\Vert w(F)(t)\Vert \le \int^t_0 \Vert K(t-s) F(s)\Vert ds
\le C\int^t_0 (t-s)^{\alpha-1}\Vert F(s)\Vert ds, \quad 0<t<T.
$$
The Young inequality on the convolution implies
\begin{align*}
& \Vert w(F)\Vert_{L^q(0,T;X)}
= \left( \int^T_0 \Vert w(F)(t)\Vert^q dt \right)^{\frac{1}{q}}\\
\le& C\Vert t^{\alpha-1}\Vert_{L^1(0,T)}
\left( \int^T_0 \Vert F(s) \Vert^q dt \right)^{\frac{1}{q}}
= \frac{CT^{\alpha}}{\alpha}\Vert F\Vert_{L^q(0,T;X)}.
\end{align*}
Thus the proof of Lemma 6.4 is complete.
$\blacksquare$

We proceed to 
\\
{\bf Completion of Proof of Theorem 6.1.}
\\
We are given $F\in L^q(0,T;X)$ arbitrarily.  Then we can find a
sequence $F_n \in C^{\infty}_0(\OOO\times (0,T))$, $n\in \N$ such that 
$\lim_{n\to\infty} \Vert F_n - F\Vert_{L^q(0,T;X)} = 0$.
Theorem 4.1 (1) yields $w(F_n) \in \DDD(\pppa) \cap L^q(0,T;\DDD(A))$.
Hence, we can readily verify 
$$
\langle w(F_n), \, ((\pppa)' - A')\psi \rangle
= \langle F_n,\, \psi \rangle \quad \mbox{for all
$\psi \in \Psi$ and $n\in \N$.}
$$
On the other hand, Lemma 6.4 implies
$$
w(F_n) \longrightarrow w(F) \quad \mbox{in $L^q(0,T;X)$ as 
$n\to \infty$}.            \eqno{(6.7)}
$$
In terms of (6.7), we can let $n \to \infty$, so that 
$$
\langle w(F), \, ((\pppa)' - A')\psi \rangle
= \langle F,\, \psi \rangle \quad \mbox{for all
$\psi \in \Psi$.}
$$
This means that $w(F)$ is a weak solution to (6.1), so that the 
existence of a weak solution $u$ is proved.
Thus the proof of Theorem 6.1 is complete.
$\blacksquare$
\\
\vspace{0.2cm}

{\bf 6.2. Smoother solution}

We consider smoother solutions, assuming that $a=0$.
More precisely, we discuss 
$$
\pppa u(t) = Au(t) + F(t), \quad 0<t<T.      \eqno{(6.8)}
$$
Similarly to the case $\alpha=1$ (e.g., \cite{Pa}, \cite{Ta}), as a class of 
$F$, we introduce a function space of H\"older continuous functions
$C^{\sigma}([0,T];X)$ with $\sigma \in (0,1)$ to specify the 
regularity of the solution to (6.8).  More precisely, we set 
$$
C^{\sigma}([0,T];X):= \left\{ u\in C([0,T];X);\, 
\sup_{0\le t \le s\le T} \frac{\Vert u(s) - u(t)\Vert}{\vert s-t\vert
^{\sigma}} < \infty\right\}.
$$
We define the norm by 
$$
\Vert u\Vert_{C^{\sigma}([0,T];X)}
:= \Vert u\Vert_{C([0,T];X)} 
+ \sup_{0\le t \le s\le T} \frac{\Vert u(s) - u(t)\Vert}{\vert s-t\vert
^{\sigma}},
$$
and we see that $C^{\sigma}([0,T];X)$ is a Banach space (e.g., 
\cite{Ya}).

We are ready to state the main result in this subsection.
\\
{\bf Theorem 6.2.}
\\
{\it
We assume 
$$
\sigma < 1 - \alpha.                        \eqno{(6.9)}
$$
Let $F \in C^{\sigma}([0,T];X)$ and $F(0) = 0$.  Then, 
$$
u(t) = \int^t_0 K(t-s)F(s) ds
$$
satisfies $Au, \pppa u \in C^{\sigma}([0,T];X)$ and (6.8) for  
$0<t<T$.  Moreover, there exists a constant $C>0$ such that 
$$
\Vert Au\Vert_{C^{\sigma}([0,T];X)}
+ \Vert \pppa u\Vert_{C^{\sigma}([0,T];X)}
\le C\Vert F\Vert_{C^{\sigma}([0,T];X)}.
$$
}

We do not know whether the conclusion of the theorem holds true also for
$\sigma \ge 1 - \alpha$.
In the case of $\alpha=1$, it is a classical result that the corresponding 
conclusion 
is satisfied for all $0 < \sigma < 1$ (e.g., \cite{Pa}).
Moreover if $q=2$, $X$ is a Hilbert space and $A$ is self-adjoint with 
compact resolvent, then the theorem holds for each 
$0<\sigma<1$ and $0<\alpha<1$ (e.g., Theorem 2.4 (3) in \cite{SY}).

Choosing a Banach space $Z$ as space of solutions, 
in general we say that the maximum regularity 
in $Z$ for (6.8) holds if $F\in Z$, then there exists a unique solution 
$u$ and $Au, \, \pppa u \in Z$.  Thus Theorem 6.2 asserts 
that the maximum regularity holds for 
(6.8) where $Z = C^{\sigma}([0,T];X)$ with a Banach space and 
$\sigma < 1-\alpha$.

With different choices of $Z$ such as $Z= L^q(0,T;X)$, there are many 
results for $\alpha=1$ and we can here refer only to Yagi \cite{Ya} and
the references therein.
On the other hand, for $\alpha \ne 1, 2$, there are few results.
\\
\vspace{0.1cm}
\\
{\bf Proof of Theorem 6.2.}
\\
{\bf First Step: estimation of $Au(t+h) - Au(t)$.}
\\
The proof is based on a similar idea to the case $\alpha=1$
(e.g., Lemma 3.4 (p.113) and Theorem 3.5 (p.114) in \cite{Pa}).

Without loss of generality, we can assume that $0\le t\le T$ and $h>0$.
Setting  
$$
v_1(t) = \int^t_0 K(t-\xi)(F(\xi) - F(t)) d\xi, \quad
v_2(t) = \int^t_0 K(t-\xi) F(t) d\xi,
$$
we have 
$$
u = v_1 + v_2 \quad \mbox{in $(0,T)$}.
$$
We will separately estimate $v_1$ and $v_2$.
Then,
\begin{align*}
& Av_1(t+h) - Av_1(t) \\
= & \int^t_0 \{ K(t+h-\xi)(F(\xi) - F(t+h))
- K(t-\xi)(F(\xi) - F(t)) \} d\xi\\
+& \int^{t+h}_t K(t+h-\xi)(F(\xi) - F(t+h)) d\xi.
\end{align*}
Moreover, we decompose 
\begin{align*}
&  \int^t_0 (K(t+h-\xi)(F(\xi) - F(t+h))
- K(t-\xi)(F(\xi) - F(t)) d\xi\\
=& \int^t_0 \{ K(t+h-\xi) \{(F(\xi) - F(t)) + (F(t) - F(t+h))\} \} d\xi 
- \int^t_0 K(t-\xi)(F(\xi) - F(t)) d\xi\\
= & \int^t_0 (K(t+h-\xi) - K(t-\xi))(F(\xi) - F(t)) d\xi
+ \int^t_0 K(t+h-\xi)(F(t) - F(t+h)) d\xi
\end{align*}
Hence,
\begin{align*}
& Av_1(t+h) - Av_1(t)
= \int^t_0 A(K(t+h-\xi) - K(t-\xi))(F(\xi) - F(t)) d\xi \\
+& \int^t_0 AK(t+h-\xi)(F(t) - F(t+h)) d\xi
+ \int^{t+h}_t AK(t+h-\xi)(F(\xi) - F(t+h)) d\xi \\
=: & S_1 + S_2 + S_3.
\end{align*}
We set $M: = \Vert F\Vert_{C^{\sigma}([0,T];X)}$ and so 
$\Vert F(\xi) - F(t)\Vert \le M\vert s-t\vert^{\sigma}$.

In terms of Lemma 3.4, we estimate
\begin{align*}
& \Vert S_1\Vert \le \int^t_0 \Vert AK(t+h-\xi) - AK(t-\xi)\Vert 
\Vert F(\xi) - F(t)\Vert d\xi\\
\le& C\int^t_0 (t-\xi)^{-\alpha}((t-\xi)^{\alpha-1} 
- (t-\xi+h)^{\alpha-1}) M\vert t-\xi\vert^{\sigma} d\xi  \\
= & CM\int^t_0 \eta^{\sigma-\alpha}(\eta^{\alpha-1} - (\eta+h)^{\alpha-1}) 
d\eta.
\end{align*}
In the last equality we use the change of the variables: $\eta = t - \xi$.
Since $(a+b)^{1-\alpha} \le a^{1-\alpha} + b^{1-\alpha}$ for 
$a,b, \ge 0$, we can continue as 
$$
 \int^t_0 \eta^{\sigma-\alpha}(\eta^{\alpha-1} - (\eta+h)^{\alpha-1}) d\eta
= \int^t_0 \eta^{\sigma-1}\left( 1 - \left( \frac{\eta}{\eta + h}\right)
^{1-\alpha}\right) d\eta 
$$
$$
\le \int^t_0 \eta^{\sigma-1}\left( \frac{h}{\eta+h}\right)
^{1-\alpha} d\eta
= h^{\sigma}\int^{\frac{t}{h}}_0 s^{\sigma-1}\left( \frac{1}{1+s}\right)
^{1-\alpha} ds.                                    \eqno{(6.10)}
$$
Here we change the variables: $s = \frac{\eta}{h}$.
Consequently,
$$
 \int^t_0 \eta^{\sigma-\alpha}(\eta^{\alpha-1} - (\eta+h)^{\alpha-1}) d\eta
\le h^{\sigma}\int^{\infty}_0 s^{\sigma-1}(1+s)^{\alpha-1}ds
< \infty
$$
by (6.9) and $\sigma > 0$.
Hence, we reach $\Vert S_1\Vert \le CMh^{\sigma}$.

Next, by (4.3), we have
\begin{align*}
& S_2 = \left( \int^t_0 A\left( \frac{d}{d\eta}J^{\alpha}G\right)(t+h-\xi) 
d\xi \right)(F(t) - F(t+h))\\
=& -\left( \int^t_0 A \frac{d}{d\xi}(J^{\alpha}G(t+h-\xi)) d\xi\right)
(F(t) - F(t+h))  \\
= & (AJ^{\alpha}G(t+h) - AJ^{\alpha}G(h))(F(t) - F(t+h)),
\end{align*}
and so Lemma 3.4 (ii) yields
$$
\Vert S_2\Vert \le C\Vert F(t+h) - F(t)\Vert 
\le CMh^{\sigma}.
$$
Finally, by (4.9) with $\beta = 1$, using $\Vert F(t+h) - F(\xi)\Vert
\le CM\vert t+h-\xi\vert^{\sigma}$, we can estimate
\begin{align*}
& \Vert S_3\Vert \le C\int^{t+h}_t (t+h-\xi)^{-1}
\Vert F(t+h) - F(\xi)\Vert d\xi\\
\le& CM\int^{t+h}_t (t+h-\xi)^{-1} (t+h-\xi)^{\sigma} d\xi
= \frac{CM}{\sigma} h^{\sigma}.
\end{align*}
Summing up, we reach 
$$
\Vert Av_1(t+h) - Av_1(t)\Vert 
\le C\Vert F\Vert_{C^{\sigma}([0,T];X)}h^{\sigma}
                                             \eqno{(6.11)}
$$
for all $F \in C^{\sigma}([0,T];X)$.
\\
{\bf Second Step.}
\\
We prove
$$
A\left( \int^t_0 K(\xi) d\xi\right)a = G(t)a - a
\quad \mbox{for all $a \in \DDD(A)$}.         \eqno{(6.12)}
$$
\\
{\bf Proof of (6.12).}
\\
By (4.3) and (2.10), we have 
\begin{align*}
& A\left( \int^t_0 K(\xi) d\xi\right)a 
= A\int^t_0 \frac{d}{d\xi}J^{\alpha}G(\xi)a \, d\xi\\
=& \int^t_0 \frac{d}{d\xi}J^{\alpha}AG(\xi)a \, d\xi
= \int^t_0 \frac{d}{d\xi}J^{\alpha}\ppp_{\xi}^{\alpha}(G(\xi)a-a)\, d\xi.
\end{align*}

On the other hand, the definition of $\pppa$ means that 
$g\in W_{\alpha,q}(0,T;X)$ 
implies that $g=J^{\alpha}w$ with some $w \in L^q(0,T;X)$ and 
$w(\xi) = \ppp_{\xi}^{\alpha}g(\xi)$.
Therefore, 
$$
J_{\xi}^{\alpha}\ppp_{\xi}^{\alpha} g(\xi) = J_{\xi}^{\alpha}w(\xi)
= g(\xi) \quad 
\mbox{for $g \in W_{\alpha,q}(0,T;X)$.}
$$
Consequently,
$$
A\left( \int^t_0 K(\xi) d\xi \right)a 
= \int^t_0 \frac{d}{d\xi}(G(\xi)a-a)\, d\xi
= G(t)a - a.
$$
Thus the proof of (6.12) is completed.
$\blacksquare$

Foe the moment we assume $F \in C^{\sigma}([0,T];\DDD(A))$.
Then, (6.12) implies 
$$
Av_2(t) = A\int^t_0 K(t-\xi)F(t) d\xi 
= A\left( \int^t_0 K(\xi) d\xi \right) F(t) = G(t)F(t) - F(t).
$$
Still we can assume that $0<t<T$ and $h>0$ is sufficiently small.
Then,
$$
Av_2(t+h) - Av_2(t) 
= G(t+h)F(t+h) - G(t)F(t) - (F(t+h)-F(t)).    \eqno{(6.13)}
$$
By Lemma 2.1 (ii) with $\beta = 0$, we obtain
\begin{align*}
& \Vert G(t+h)F(t+h) - G(t)F(t)\Vert
= \Vert G(t+h)(F(t+h)-F(t)) + (G(t+h) -G(t))F(t)\Vert\\
\le& \Vert G(t+h)\Vert \Vert F(t+h)-F(t)\Vert 
+ \Vert (G(t+h) -G(t))F(t)\Vert
\end{align*}
$$
\le CMh^{\sigma} + \Vert (G(t+h) -G(t))F(t)\Vert.
                                                      \eqno{(6.14)}
$$
On the other hand, since $F(0) = 0$, the application of Lemma 3.1 
with $\beta = 0$ yields
\begin{align*}
& \Vert (G(t+h) -G(t))F(t)\Vert
= \left\Vert \left( \int^{t+h}_t \frac{dG}{d\xi}(\xi) d\xi \right)
(F(t) - F(0)) \right\Vert
\le \int^{t+h}_t \xi^{-1}\Vert F(t) - F(0)\Vert d\xi\\
\le& M\int^{t+h}_t \xi^{-1}t^{\sigma} d\xi
\le M\int^{t+h}_t \xi^{-1}\xi^{\sigma} d\xi 
= \frac{M}{\sigma}((t+h)^{\sigma} - t^{\sigma}) \le \frac{M}{\sigma}h^{\sigma},
\end{align*}
where we used $(t+h)^{\sigma} \le t^{\sigma} + h^{\sigma}$.
Therefore, (6.13) and (6.14) imply 
$$
\Vert Av_2(t+h) - Av_2(t)\Vert 
\le C\Vert F\Vert_{C^{\sigma}([0,T];X)} h^{\sigma},
$$
that is,
$$
\Vert Au(t+h) - Au(t)\Vert 
\le C\Vert F\Vert_{C^{\sigma}([0,T];X)} h^{\sigma}.   \eqno{(6.15)}
$$
Moreover we can easily prove $\Vert Au(t)\Vert \le
C\Vert F\Vert_{C^{\sigma}([0,T];X)}$.
Hence,
$$
\Vert Au\Vert_{C^{\sigma}([0,T];X)} + \Vert \pppa u\Vert_{C^{\sigma}([0,T];X)}
\le C\Vert F\Vert_{C^{\sigma}([0,T];X)}              \eqno{(6.16)}
$$
for all $F \in C^{\sigma}([0,T];\DDD(A))$.
\\
Thus the proof of Theorem 6.2 is complete if $F \in C^{\sigma}([0,T];\DDD(A))$.
\\
{\bf Third Step.}
\\
We will complete the proof of Theorem 6.2 for $F \in C^{\sigma}([0,T];X)$.
We can choose a sequence $F_n\in C^{\sigma}([0,T];\DDD(A))$, $n\in \N$ 
such that $\Vert F_n - F\Vert_{C^{\sigma}([0,T];X)} \rrrr 0$ as 
$n\to \infty$.  We set $u(F)(t) = \int^t_0 K(t-s)F(s) ds$.
Then, we can see
$$
\pppa u(F_n)(t) = Au(F_n)(t) + F_n(t), \quad 0<t<T.     \eqno{(6.17)}
$$
In view of (6.16), there exists $\www{u} \in C^{\sigma}([0,T];X)$
such that $\lim_{n\to \infty} \Vert u(F_n) - \www{u}\Vert
_{C^{\sigma}([0,T];X)} = 0$ and 
$Au(F_n)$ and $\pppa u(F_n)$ converge in $C^{\sigma}([0,T];X)$.
Since $\pppa: \DDD(\pppa) \subset L^q(0,T;X) \rrrr L^q(0,T;X)$ and
$A: L^q(0,T;\DDD(A)) \subset L^q(0,T;X) \rrrr L^q(0,T;X)$ are closed 
operators, letting $n\to \infty$, we obtain 
$\pppa \www{u}(t) = A\www{u}(t) + F(t)$ in
$(0,T)$.  Moreover by (4.9) with $\beta = 0$, we have
\begin{align*}
& \left\Vert \int^t_0 K(t-s)F_n(s) ds - \int^t_0 K(t-s)F(s)ds \right\Vert
\le C\int^t_0 (t-s)^{\alpha-1}\Vert F_n(s) - F(s)\Vert ds\\
\le &C\left( \int^t_0 (t-s)^{\alpha-1} ds\right) \Vert F_n - F\Vert
_{L^{\infty}(0,T;X)} 
\le CT^{\alpha}\Vert F_n - F\Vert_{L^{\infty}(0,T;X)}
\, \rrrr \, 0
\end{align*}
for $0\le t \le T$.
Therefore, $\www{u}(t) = \int^t_0 K(t-s)F(s) ds$ for $0<t<T$ and 
$$
\Vert A\www{u}\Vert_{C^{\sigma}([0,T];X)}
+ \Vert \pppa \www{u}\Vert_{C^{\sigma}([0,T];X)} 
\le C\Vert F\Vert_{C^{\sigma}([0,T];X)}.
$$
Thus the proof of Theorem 6.2 is complete for all $F \in 
C^{\sigma}([0,T];X)$.
$\blacksquare$
\section{Application to an inverse problem}

There are various inverse problems for time-fractional partial 
differential equations and we here discuss one inverse problem
within our framework in order to demonstrate the wide applicability 
of our method.  There are many works on inverse problems with 
the framework of $L^2(\OOO)$ or Hilbert spaces, but there are very few 
publications within non-Hilbert spaces.

One of the main method for inverse problems is the Laplace transform.
Our construction of solution to an initial boundary value problem
is based on the Laplace transform, so that the application of 
Laplace transforms for inverse problems is quite natural.

We consider the case where $A$ is an elliptic operator given by 
(1.7) and (1.8), where we assume that the $b_k = 0$ for $1 \le k \le d$
and $a_{k\ell}, c \in C^{\infty}(\ooo{\OOO})$.
We assume that $1<p<\infty$.
For 
$$
\pppa (u(t) - a) = Au(t), \quad t>0,         \eqno{(7.1)}
$$
we consider
\\
{\bf Inverse problem of determining initial value.}
\\
{\it 
Let $a \in L^p(\OOO)$ and let $\omega \subset \OOO$ be an arbitrarily 
given subdomain.  Then, determine $a(x)$ for $x \in \OOO$ by 
$u\vert_{\omega \times (0,T)}$.
}
\\

We state our main results in this section.
\\
{\bf Theorem 7.1.}
\\
{\it
Let $T>0$ be arbitrarily chosen.  Then, $u\vert_{\omega \times (0,T)} = 0$ 
implies that $a=0$ in $\OOO$.
}
\\

The theorem means that if 
$\pppa (u-a) = Au$ and $\pppa (\www{u} - \www{a}) = A\www{u}$ in 
$\OOO\times (0,T)$, and $u = \www{u}$ in $\omega \times (0,T)$, then 
$a = \www{a}$ in $\OOO$.
\\

In fact, we can prove a stronger uniqueness than Theorem 7.1.
We further introduce notations.  For $1<p<\infty$, by $A_p$ we denote the 
elliptic operator defined by (1.7) with the domain 
$\DDD(A_p) = \{ v\in W^{2,p}(\OOO);\, w\vert_{\ppp\OOO} = 0\}$.
Here and henceforth we write $A_p$, $\AAAAT$ in order to 
specify that the working space is $L^p(\OOO)$.  
We introduce another elliptic operator similarly: let 
$$
\www{\AAA}v(x) = \sum_{k,\ell=1}^d \ppp_k(\www{a}_{k\ell}(x)\ppp_{\ell}v(x))
+ \www{c}(x)v(x), \quad x\in \OOO.
                                     \eqno{(7.2)}
$$
By $\www{A}_p$ we define an operator $\AAA$ in $L^p(\OOO)$ attached with 
the domain $\DDD(\www{A}_p) = \{ v\in W^{2,p}(\OOO);\, 
v\vert_{\ppp\OOO} = 0\}$.  Here $\www{a_{k\ell}} = \www{a_{\ell k}}$,
$c \in C^{\infty}(\ooo{\OOO})$ and (1.8) is satisfied 
by $\www{a_{k\ell}}$ and $\www{c}$.  Then
\\
{\bf Theorem 7.2.}
\\
{\it
We assume that $A_p$ and $\www{A}_p$ are commutative, that is,
$$
A_p\www{A}_p = \www{A}_pA_p, \quad
\{ v\in \DDD(\www{A}_p);\, \www{A}_pv \in \DDD(A_p)\}
= \{ v\in \DDD(A_p);\, A_pv \in \DDD(\www{A}_p)\}.              \eqno{(7.3)}
$$
Let 
$$
\left\{ \begin{array}{rl}
& \pppa (u(t) - a) = A_pu(t), \quad 0<t<T, \\
& \pppa (\www{u}(t) - \www{a}) = \www{A}_p\www{u}(t), \quad 0<t<T.
\end{array}\right.
                                                \eqno{(7.4)}
$$
If $u = \www{u}$ in $\omega \times (0,T)$, then $a = \www{a}$ in $\OOO$.
}


As $\AAA$ and $\www{\AAA}$ satisfying (7.3), we can mention that 
all the coefficients of $\AAA$ and $\www{A}$ are constants.
Therefore,
\\
{\bf Corollary 7.1.}
\\
{\it
Let $A_p$ and $\www{A}_p$ be defined through (1.7) and (7.2).
We assume that $a_{k\ell}, c, \www{a_{k\ell}}, \www{c}$ are all constants.
Then the same assertion as Theorem 7.2 holds true.
}

Theorem 7.2 means that we can prove the uniqueness of initial value
although we do not know the elliptic operator $\AAA$ and $\www{\AAA}$.
In other words, the initial value is much more distinguishable than 
any other quantities related to the equations.

Theorem 7.1 directly follows by choosing $A = \www{A}_p$ in Theorem 7.2.
Hence, it suffices to prove Theorem 7.2.
\\
{\bf Proof of Theorem 7.2.}
\\
In terms of the time analyticity of $u(t)$ and $\www{u}(t)$, we can obtain 
that (7.4) holds for $t>0$.
Therefore, by Theorem 4.1 (2), we have
$$
(Lu)(\la) = \la^{\alpha-1}(\la^{\alpha}-A_p)^{-1}a
$$
and
$$
(L\www{u})(\la) = \la^{\alpha-1}(\la^{\alpha}-\AAAAT)^{-1}a
\quad \mbox{in $L^p(\OOO)$ if $\mbox{Re}\, \la > \la_0$ and
$\la^{\alpha} \in \rho(A_p) \cap \rho(\AAAAT)$.}
$$
Here $\la_0>0$ is some constant.
Hence, by means of the time-analyticity, from 
$u=\www{u}$ in $\omega \times (0,T)$, we have 
$u=\www{u}$ in $\omega \times (0,\infty)$, that is,
$$
\la^{\alpha-1}(\la^{\alpha}-A_p)^{-1}a(x) 
= \la^{\alpha-1}(\la^{\alpha}-\AAAAT)^{-1}\www{a}(x) \quad
\mbox{for $x \in \omega$},
$$
if $\mbox{Re}\, \la > \la_0$, $\la^{\alpha} \in \rho(A_p)\cap
\rho(\AAAAT)$, that is,
$$
(\la^{\alpha}-A_p)^{-1}a(x) = (\la^{\alpha}-\AAAAT)^{-1}\www{a}(x) \quad
\mbox{for $x \in \omega$}
$$
if $\mbox{Re}\, \la > \la_0$, $\la^{\alpha} \in \rho(A_p) 
\cap \rho(\AAAAT)$.
Setting $z = \la^{\alpha}$, we have 
$$
(z-A_p)^{-1}a(x) = (z-\AAAAT)^{-1}\www{a}(x) \quad
\mbox{for $z \in \mbox{Re}\, z > \la_0^{\alpha}$ and 
$z \in \rho(A_p) \cap \rho(\AAAAT)$ and $x \in \omega$.}
$$
Since $(z-A_p)^{-1}a$ and $(z-\AAAAT)^{-1}a$ are holomorphic in 
$\rho(A_p)$ and $\rho(\AAAAT)$ respectively, we obtain
$$
(z-A_p)^{-1}a(x) = (z-\AAAAT)^{-1}\www{a}(x) \quad
\mbox{for $x \in \omega$ and $z \in \C \setminus (\sigma(A_p) \cup
\sigma(\AAAAT))$.}
                                             \eqno{(7.5)}
$$
It is known that the spectra $\sigma(A_p) := \C \setminus \rho(A_p)$ and 
$\sigma(\AAAAT) := \C \setminus \rho(\AAAAT)$ are both countably infinite sets
in $\R$.
Moreover $\sigma(A_p)$ and $\sigma(\AAAAT)$ are independent of 
$p\in (1,\infty)$, and especially
$$
\sigma(A_p) = \sigma(A_2), \quad \sigma(\AAAAT) = \sigma(\www{A}_2)
\quad \mbox{for $1<p<\infty$}.                               \eqno{(7.6)}
$$
\\
{\bf Proof of (7.6).}
\\
Since all the coefficients of $A_p$ are in $C^{\infty}(\ooo{\OOO})$, in 
terms of (4.15), we can prove $\DDD(A_p^m) \subset W^{2m,p}(\OOO)$ for each 
$m\in \N$.  Hence,
$$
\DDD(A_p^m) \subset \DDD(A_2) \quad \mbox{with some $m\in \N \cup \{0\}$}.
                                                \eqno{(7.7)}
$$
Indeed, for $1<p<2$, we choose large $m\in \N$, so that 
the Sobolev embedding yields (7.7). If $p \ge 2$, then (7.7) already 
holds with $m=0$.

We note that $\la \in \sigma(A_p)$ if and only if there exists
$\va \in \DDD(A_p)$, $\ne 0$ such that $A_p\va = \la \va$.
Then, $A_p^m \va = \la^m\va \in \DDD(A_p)$ for each $m\in \N$, which 
implies that $\va \in \DDD(A_p^m)$.  Therefore, if 
$\va \in \DDD(A_p)$ is an eigenfunction of $A_p$ for $\la$, then 
(7.7) yields $\va \in \DDD(A_2)$, that is, $\sigma(A_p) \subset 
\sigma(A_2)$.  The inverse inclusion can be proved similarly.
$\blacksquare$
\\

We define the eigenprojections:
$$
\PPPPO := \frac{1}{2\pi i} \int_{C_{\la}} (z-A_p)^{-1} dz, \quad 
\PPPPT := \frac{1}{2\pi i} \int_{C_{\la}} (z-\AAAAT)^{-1} dz
$$
for $\la \in \sigma(A_p) \cup \sigma(\AAAAT)$, where we specify
the index $p$.
Here and henceforth, $C_{\la}$ denotes a circle centered at $\la$
with sufficiently small radius which does not surround any 
points from $\sigma(A_p) \cup \sigma(\AAAAT) \setminus \{ \la\}$.

Let $\la \in \sigma(A_p)$.  First we consider the case
$\la \in \sigma(A_p) \setminus \sigma(\AAAAT)$.
Taking the integral $\int_{C_{\la}} \cdots dz$ of (7.5), 
we obtain $\PPPPO a(x) = 0$ for $x \in \omega$.
Since $\PPPPO a$ belongs to the eigenspace of $A_p$ for $\la$, we have
$(A-\la)\PPPPO a = 0$ in $\OOO$.  By $P_{\la}^{(p)}a = 0$ 
in $\omega$, we apply the unique continuation for 
an elliptic operator $A - \la$ in $\OOO$, so that $\PPPPO a = 0$ in $\OOO$.
Similarly, if $\la \in \sigma(\AAAAT) \setminus \sigma(A_p)$, 
then $\PPPPT \www{a} = 0$ in $\OOO$.

Second we consider the case $\la \in \sigma(A_p) \cap \sigma(\AAAAT)$.
Taking the integral along $C_{\la}$ of (7.5), we can reach 
$\PPPPO a = \PPPPT \www{a}$ in $\omega$.
We set $y_{\la}:= \PPPPO a - \PPPPT \www{a}$.
Then, 
$$
y_{\la} = 0 \quad \mbox{in $\omega$}.              \eqno{(7.8)}
$$
By the sufficient smoothness of the coefficients $a_{k\ell}, c, 
\www{a}_{k\ell}, \www{c}$ on $\ooo{\OOO}$, we see that the solution 
$\PPPPO a$ to $(A-\la)\PPPPO a = 0$ in $\OOO$ and 
$\PPPPO a = 0$ on $\ppp\OOO$, is in $W^{4,p}(\OOO)$, and 
similarly $\PPPPT \www{a} \in W^{4,p}(\OOO)$.  Therefore, 
$y_{\la} \in W^{4,p}(\OOO)$.  Hence, in view of (7.3), we have
\begin{align*}
& (\AAAAT - \la)(A_p-\la)y_{\la} 
= (\AAAAT - \la) (A_p - \la)\PPPPO a - (\AAAAT - \la)(A_p-\la)
\PPPPT\www{a}\\
=& -(A_p - \la)(\AAAAT-\la)\PPPPT \www{a}
=0 \quad \mbox{in $\OOO$}.
\end{align*}
Setting $z_\la:= (A_p-\la)y_{\la}$ in $\OOO$, by (7.8) we have also 
$z_{\la}=0$ in $\omega$.
Thus $z_{\la}\in W^{2,p}(\OOO)$ satisfies 
$(\AAAAT - \la)z_{\la}=0$ in $\OOO$.  The unique continuation for an elliptic 
operator $\AAAAT - \la$ in $\OOO$, yields $z_{\la} 
= (A_p-\la)y_\la=0$ in $\OOO$.
Since $y_{\la} = 0$ in $\omega$ by (7.8), again the 
unique continuation to $A_p- \la$ in $\OOO$,  yields 
$y_{\la} = 0$ in $\OOO$, that is,
$$
\PPPPO a = \PPPPT \www{a} \quad \mbox{in $\OOO$ for each 
$\la \in \sigma(A_p) \cap \sigma(\AAAAT)$}. 
$$

Thus we reach
$$
\left\{ \begin{array}{rl}
& \PPPPO a = 0 \quad \mbox{in $\OOO$ if $\la \in 
\sigma(A_p) \setminus \sigma(\AAAAT)$}, \cr\\
& \PPPPO a = \PPPPT \www{a} \quad \mbox{in $\OOO$ if $\la \in 
\sigma(A_p) \cap \sigma(\AAAAT)$}, \cr\\
& \PPPPT \www{a} = 0 \quad \mbox{in $\OOO$ if $\la \in 
\sigma(\AAAAT) \setminus \sigma(A_p)$}.
\end{array}\right.
                                         \eqno{(7.9)}
$$

In (7.9), we lift up the regularity of $a, \www{a}$ by a similar 
idea used in the proof of Proposition 4.4.
In view of $A_p\AAAAT = \AAAAT A_p$, we can prove 
$$
A_p^m(z-\AAAAT) = (z-\AAAAT)A_p^m \quad \mbox{in $\DDD(A_p^{m+1})$}
$$
for each $m\in \N$.  Hence, $(A_p^m(z-\AAAAT))^{-1} 
= ((z-\AAAAT)A_p^m)^{-1}$ for $z \in \rho(\AAAAT)$, that is,
$$
(z-\AAAAT)^{-1} (A_p^m)^{-1} = (A_p^m)^{-1}(z-\AAAAT)^{-1} \quad
\mbox{for $z\in \rho(\AAAAT)$.}
$$
Consequently, taking the integral along $C_{\la}$ and recalling the 
definition of $\PPPPT$, we can obtain
$$
(A_p^m)^{-1}\PPPPT = \PPPPT (A_p^m)^{-1} \quad 
\mbox{for $\la \in \sigma(\AAAAT)$, $m\in \N$ and $1<p<\infty$}.
                                                      \eqno{(7.10)}
$$
Applying (7.10) in (7.9) and $A_p\PPPPO = \PPPPO A_p$ 
(e.g., Kato \cite{Ka}) and choosing $m\in \N$ sufficiently large
such that $b:= (A_p^m)^{-1}a, \www{b} := (A_p^m)^{-1}\www{a} \in 
L^p(\OOO) \cap L^2(\OOO)$, by (7.6) and (7.9) we see 
$$
\left\{ \begin{array}{rl}
& \PPPPO b = 0 \quad \mbox{in $\OOO$ if $\la \in 
\sigma(A_2) \setminus \sigma(\www{A}_2)$}, \cr\\
& \PPPPO b = \PPPPT \www{b} \quad \mbox{in $\OOO$ if $\la \in 
\sigma(A_2) \cap \sigma(\www{A}_2)$}, \cr\\
& \PPPPT \www{b} = 0 \quad \mbox{in $\OOO$ if $\la \in 
\sigma(\www{A}_2) \setminus \sigma(A_2)$}.
\end{array}\right.
                                         \eqno{(7.11)}
$$
Moreover, since $A_p \subset A_{p'}$ and $\AAAAT \subset \www{A}_{p'}$ for 
$p' < p$, noting (7.6), by the definition of $\PPPPO$ and $\PPPPT$
we know that $\PPPPO \subset P_{\la}^{(p')}$ for $\la \in \sigma(A_2)$ and
$\PPPPT \subset \www{P}_{\la}^{(p')}$ for $\la \in \sigma(\www{A}_2)$.

Therefore, by $b, \www{b} \in L^p(\OOO) \cap L^2(\OOO)$, 
we can readily verify that 
$\PPPPO b = P^{(2)}_{\la} b$ and $\PPPPT \www{b} = \www{P}_{\la}^{(2)}
\www{b}$.  Consequently, (7.11) implies 
$$
\left\{ \begin{array}{rl}
& P_{\la}^{(2)}b = 0 \quad \mbox{in $\OOO$ if $\la \in 
\sigma(A_2) \setminus \sigma(\www{A}_2)$}, \cr\\
& P_{\la}^{(2)} b = \www{P}_{\la}^{(2)}\www{b} \quad \mbox{in $\OOO$ if $\la 
\in \sigma(A_2) \cap \sigma(\www{A}_2)$}, \cr\\
& \www{P}_{\la}^{(2)}\www{b} = 0 \quad \mbox{in $\OOO$ if $\la \in 
\sigma(\www{A}_2) \setminus \sigma(A_2)$}.
\end{array}\right.
                                         \eqno{(7.12)}
$$
Now we reduced relation (7.9) in $L^p(\OOO)$ to (7.12) in $L^2(\OOO)$.

Since we can construct orthonormal bases in $L^2(\OOO)$ composed of 
the eigenfunctions of $A_2$ and the ones of $\www{A}_2$ respectively, we have
$b = \sum_{\la\in \sigma(A_2)} P_{\la}^{(2)}b$ and
$\www{b} = \sum_{\la\in \sigma(\www{A}_2)} \www{P}_{\la}^{(2)}\www{b}$ in 
$L^2(\OOO)$.  Therefore, (7.12) yields
\begin{align*}
& b = \left( \sum_{\la\in \sigma(A_2) \setminus \sigma(\www{A}_2)}
 + \sum_{\la\in \sigma(A_2) \cap \sigma(\www{A}_2)} \right)
P_{\la}^{(2)}b
= \sum_{\la\in \sigma(A_2) \cap \sigma(\www{A}_2)} P_{\la}^{(2)}b\\
=& \sum_{\la\in \sigma(\www{A}_2) \cap \sigma(A_2)} \www{P}^{(2)}_{\la}\www{b}
= \left( \sum_{\la\in \sigma(\www{A}_2) \setminus \sigma(A_2)}
 + \sum_{\la\in \sigma(\www{A}_2) \cap \sigma(A_2)} \right)
\www{P}_{\la}^{(2)}\www{b} = \www{b} \quad \mbox{in $\OOO$}.
\end{align*}
Since $(A^m_2)^{-1}$ is injective, we reach $a = \www{a}$ in $\OOO$.
Thus the proof of Theorem 7.2 is complete.
$\blacksquare$

\section{Concluding remarks}

We consider an evolution equation 
$$
\pppa (u(t) - a) = Au(t), \quad 0<t<T             \eqno{(8.1)}
$$
in Banach space $X$ with time-fractional differential operator $\pppa$ 
of the order 
$\alpha \in (0,1)$, and $a \in X$ corresponds to an initial value.
If we choose $X=L^p(\OOO)$ with bounded domain $\OOO \subset \R^d$ and 
suitable elliptic operator $A$ attached with boundary condition, 
then (8.1) describes an initial boundary 
value problem for a time-fractional diffusion equation.
 
Most of the existing works are concerned with Hilbert space $X$, and 
our purpose is to provide unified treatments which are widely applicable 
to non-Hilbert space $X$.

The first subject is the well-posedness for initial value prolems (8.1).
Our method is based on the vector-valued Laplace transform, and is 
a modification of a classical construction of the analytic semigroup.
The construction derives solution formula which is feasible and convenient.

We intend to totally transfer the operator theoretic approach for 
the parabolic case (i.e., $\alpha=1$) to 
time-fractional evolution equations. 
In this article, we established at least the core part of the 
expected theory for time-fractional evolution equations.  

Such a fundamental theory for the well-posedness should be 
feasible and applicable, and we demonstrate its applicability 
in Sections 5 - 7.
The treated topics in these sections can be much more comprehensively 
studied, but in this article, we are restricted to partial studies in order to 
show the essence of the applications of the fundamental theory 
in Sections 1 - 4.   
\\
\vspace{0.2cm}
 
{\bf Acknowledgements.}
Masahiro Yamamoto was supported by 
Grant-in-Aid for Scientific Research (A) 20H00117 
and Grant-in-Aid for Challenging Research (Pioneering) 21K18142 of 
Japan Society for the Promotion of Science.
Most of this work has been carried out when Masahiro Yamamoto was 
a visiting professor at Sapienza University of Rome 
in October - November 2024.

\end{document}